\documentclass[14pt,oneside,a4paper,openany]{article}

\usepackage{qxy-notation}

\usepackage[T1]{fontenc}    % use 8-bit T1 fonts
\usepackage{hyperref}       % hyperlinks
\usepackage{url}            % simple URL typesetting
\usepackage{booktabs}       % professional-quality tables
\usepackage{amsfonts}       % blackboard math symbols
\usepackage{nicefrac}       % compact symbols for 1/2, etc.
\usepackage{microtype}      % microtypography
\usepackage{authblk}
\usepackage{lettrine}
\usepackage{csquotes}

\usepackage{subfigure}
\usepackage{color, xcolor}
\usepackage{cleveref}
\crefformat{equation}{(#2#1#3)}
\crefrangeformat{equation}{(#3#1#4) to~(#5#2#6)}

\usepackage[english]{babel}

\RequirePackage{mathtools}

\usepackage{bm}
\usepackage{natbib}
\newcounter{Lem}[section]
\renewcommand{\theLem}{\thesection.\arabic{Lem}}
\newcounter{Exam}[section]
\renewcommand{\theExam}{\thesection.\arabic{Exam}}
\newcounter{Thm}[section]
\renewcommand{\theThm}{\thesection.\arabic{Thm}}
\newcounter{Prop}[section]
\renewcommand{\theProp}{\thesection.\arabic{Prop}}
\title{On the connections between algorithmic regularization and penalization for convex losses}

\author[1]{Qian Qian \thanks{qian.216@osu.edu} }
\author[2]{Xiaoyuan Qian \thanks{xyqian@dlut.edu.cn}}
\affil[1]{Department of Statistics, Ohio State University}
\affil[2]{School of Mathematical Sciences, Dalian University of Technology}

\date{\vspace{-5ex}}

\begin{document}

\maketitle

\begin{abstract}
  In this work we establish the equivalence of algorithmic regularization and explicit convex penalization for generic convex losses.
  We introduce a geometric condition for the optimization path of a convex function, and show that if such a condition is satisfied, the optimization path of an iterative algorithm on the unregularized optimization problem can be represented as the solution path of a corresponding penalized problem.

\end{abstract}

\section{Introduction}

In statistics, estimation is often cast in terms of minimizing a loss function:
\begin{equation}\label{01001}
  \underset{x}{\arg\min}\, f(x)\,,
\end{equation}
However, direct minimization can lead to overfitting.
Instead of minimizing the loss function in \eqref{01001}, explicit penalization deals with the following optimization problem:
\begin{equation}\label{01002}
  \underset{x}{\arg \min} \ \left( f(x) + \lambda \psi(x) \right) \ ,
\end{equation}
where $\psi$ is the penalty function and $ \lambda $ is the tuning parameter.
For example, ridge regression (\cite{hoerl62}, \cite{hoerl68}), Lasso (\cite{tib96}) and elastic net (\cite{zou07}) are well-known examples of explicit regularization, with
$\psi(\theta) = { \| \theta \|}^2_2$ for ridge regression,  $ \psi(\theta) = {\|\theta \|}^2_1$ for lasso, and
$\psi(\theta) = \alpha{\|\theta \|}^2_1 + (1 - \alpha){\|\theta \|}^2_2$ where $0 < \alpha < 1$, for elastic net.
However, penalization approach requires one to solve the problem \eqref{01002} for a sequence of the tuning parameter $ \lambda $ to obtain an entire solution path,
thus yielding a considerable computational burden.
\cite{efron04} showed that the optimal solution path of Lasso is piecewise linear and proposed LARS algorithm to compute the full solution path of Lasso efficiently.
This result was extended to more generic cases by \cite{rosset07} who derived a general characterization of the properties of (loss $f$, penalty $\psi$) pairs giving piecewise linear coefficient paths that allow for efficient generation of the full regularized coefficient paths. However, this generalization holds only when the loss function $f$ is piecewise quadratic and the penalty $\psi$ is piecewise linear.
Thus the class of ($f$ and $\psi$) pairs that can be computed efficiently is limited.

Recently, there is a growing interest in the study of \emph{algorithmic regularization}: one can use an optimization algorithm (such as gradient descent) to find estimators without employing any explicit penalization.
Yet this optimization algorithm still exhibits an effect of regularization. Such regularization effects may depend on the choice of the algorithm, the loss function, the initialization and the distribution of the data. The characteristic of this algorithmic approach is that
the employed algorithm seems to perform regularization, although no explicit regularization is enforced. Therefore, in order to understand how the optimization  procedure itself affects the learned model, it is important to precisely characterize algorithmic regularization induced by different optimization techniques .

One way to study algorithmic regularization is to make connections with explicit penalization. More specifically, we can study algorithmic regularization by investigating the connection between iterates generated by optimization techniques on un-regularized objectives and minimizers of corresponding penalized objectives. These connections may help us to transfer
insights from algorithmic regularization to explicit penalization and vice versa.

\cite{friedman04} empirically observed that several methods of generalized gradient descent
are seen to produce paths that closely correspond to those induced by commonly used penalization methods. \cite{hastie2009} noted a connection between $L_2$ boosting with componentwise linear regression and Lasso.  \cite{efron04} considered the forward stagewise linear regression, which is a version of $L_2$ boosting with infinitesimally small step sizes, and show that the solutions produced by forward stagewise linear regression is equivalent to the Lasso solution path produced by varying $\lambda$. \cite{rosset04} showed that under certain conditions on the problem, the path traced by coordinate descent or boosting is similar to the regularization path of $L_1$ constrained problem. More specifically, for exponential loss and binomial log-likelihood, the boosting estimators converges to the ``$L_1$-optimal'' solution that maximizes the $L_1$ margin for separable data.
In this sense, boosting is similar to support vector machines since both methods can be viewed as regularized optimization in the predictor space.
While support vector machines solve the optimization problem exactly, boosting only solves the corresponding optimization problems approximately.

Besides $L_1$ penalization, there is also a rich literature on the connections between
early stopping of gradient descent and $L_2$ penalization.
Several works (\cite{fleming90}, \cite{santos96} and \cite{skouras94})
show that there exists a one-to-one correspondence
between the early stopping of gradient descent method on least square problems and ridge regression. Similarly, for stochastic gradient descent,
\cite{neu17} proposed a a variant of the Polyak\textendash Ruppert averaging scheme, and proved that in the context of linear least square regression,
this averaging scheme with decaying weights in a geometric fashion has the
same regularization effect, and is asymptotically equivalent to ridge regression.
  %The geometric Polyak\textendash Ruppert averaging update, given by $\tilde{\theta}_t = \frac{1}{ \sum^n_{k=0} {(1 - \eta \lambda)}^k }
   %\sum^n_{k = 0} {(1 - \eta \lambda)}^k \theta_k $ , shares a similar excess-risk
  % bound with that derived from Tikhonov regularization compared to the previous work by \cite{dieuleveut17}.
More recently, \cite{suggala2018} made connections between the optimization path of gradient descent and the corresponding $L_2$ penalization path for strongly convex training objectives.
Such a connection can also be extended to mirror descent for strongly convex loss. Moreover, a similar result also exists for unregularized logistic regression loss with separable data, which is the same situation considered in \cite{soudry17}. However, it is not known (1) if similar connections hold for  general convex losses or nonconvex losses; (2) if similar connections hold for  methods other than gradient descent, such as steepest descent, Newton's method and stochastic gradient descent.

In a general view on penalization and algorithmic regularization, one can naturally ask the following question:
under what condition there exists an equivalence between these two approaches of regularization? Or,
more practically,
is there a way that we can characterize the searching path of an iterative algorithm via
a penalization course of the loss function?
In this work we will answer the question in both necessary and sufficient aspects.
Precisely, we state a geometric condition and give the following results:

i) if %a searching path of \eqref{01001} can be characterized by a penalization course of the loss function, i.e.,
for each point at a given searching path, there is
a %positive number
$ \lambda > 0 $ such that the point is the solution of the corresponding penalization problem \eqref{01002}, then the searching path has to satisfy the geometric condition;

ii) if a discrete searching path satisfies the geometric condition, then there is a convex function $\psi$
such that for each point at the path, there is
a $ \lambda > 0 $ such that the point is the solution of \eqref{01002};

iii) if a continuous searching path satisfies the geometric condition, then for any $\varepsilon>0$, there is a convex function $\psi$ such that for each point at the searching path, there is
a $ \lambda > 0 $ such that the point is in the $\varepsilon \,$-neighborhood of the solution of \eqref{01002}.

The paper is organized as follows.
Section 2 includes notation and some assumptions used through out the paper.
In Section 3 we give a geometric characterization of the searching path that can be produced from the solutions of minimization problems of type \eqref{01002} by changing the values of the parameter $ \lambda $.
In Section 4 we show that for a discrete searching path that bears the geometric characterization, then there exists a convex function $\psi$
such that each point at the path can be obtained by solving a minimization problem of type \eqref{01002} with an appropriate $ \lambda $.
This result is extended to continuous searching paths in an approximate form in Section 5.
In Section 6 we make a short review on our results.
In addition, we leave most proofs of the lemmas in Appendix to focus our attention on the main results.

%%%%%%%%%%%%%%%%%%%%%%%%%%%%%%%%%%%%%%%%%%%%%%%%%%%%
\section{Preliminaries}
%%%%%%%%%%%%%%%%%%%%%%%%%%%%%%%%%%%%%%%%%%%%%%%%%%%%

The closure, interior, and boundary of a set $A\subset\Rn$ are denoted by
$ \mathrm{cl}\,A,\ \ \mathrm{int}\,A,$ and $\mathrm{bdry}\,A,$ respectively.
The affine hull, convex hull, and conic hull of $A$ are denoted by
$ \mathrm{aff}\,A,\ \ \mathrm{conv}\,A,$ and $\mathrm{cone}\,A,$ respectively.

The relative interior of the set $A$ is denoted by
$\mathrm{ri}\,A$.
The relative boundary of $A$ is defined as the relative complement of
$\mathrm{ri}\,A$ with respect to $\mathrm{cl}\,A $ and denoted by $ \mathrm{rbd}\,A $, i.e.
$\mathrm{rbd}\,A = (\mathrm{cl}\,A) \setminus \mathrm{ri}\,A $.

For two points $ \bm{a}$ and $\bm{b} $ in a Euclidean space,
the line segment connecting $ \bm{a}$ and $\bm{b} $ is denoted by
$ \overline{\bm{a}\bm{b}} $, i.e.
$$
  \overline{\bm{a}\bm{b}} = \{ (1-t)\bm{a}+ t\bm{b} \ |\ t\in [0,1] \}\,;
$$
and the ray starting from $ \bm{a}$ and passing through $\bm{b} $ is denoted by
$$
  \mathrm{ray}(\bm{a},\bm{b}) = \{ (1-t)\bm{a}+ t\bm{b} \ |\ t \geq 0 \}\,.
$$

Throughout this paper, we suppose that $f: \Rn\ra \exreal $ is a proper convex function.
i.e. $f$ is convex, $f(\bm{x})<\pinf $ for at least one $\bm{x}$, and $f(\bm{x})>\ninf $ for all $\bm{x}$.
The effective domain of $f$, denoted by $ \mathrm{dom}\,f $, is defined as
$$
 \mathrm{dom}\,f := \{\bm{x}\in \Rn\ |\ f(\bm{x})<\pinf\} \,.
$$
We also suppose that the set of minimizers of $f$ is nonempty,
which means $ \mathrm{dom}\,f $ is also a nonempty set.
Moreover, we assume that $f(\bm{x})=\pinf $ for all $\bm{x}\in \mathrm{bdry}\,\mathrm{dom}\,f $, which imples $ \mathrm{dom}\,f $ is open.

The symbols
$$
  \mathrm{lev}_{\leq c}f := \{ \bm{x}\in \mathrm{dom}\,f \ |\ f(\bm{x})\leq c \}
$$
and
$$
  \mathrm{lev}_{< c}f := \{ \bm{x}\in \mathrm{dom}\,f \ |\ f(\bm{x})< c \}
$$
are used to denote the lower level set and strict lower level set of a function $f$,
\rsp.

Let $\bm{x}\in \mathrm{dom}\,f$ and $ \bm{x}^{\ast} \in \partial f(\bm{x}) $.
We denoted by $ H^{+}_{f}(\bm{x}, \bm{x}^{\ast}) $,
$ H^{-}_{f}(\bm{x}, \bm{x}^{\ast}) $, and
$ H(\bm{x}, \bm{x}^{\ast}) $
the halfspaces
$$
  H^{+}(\bm{x}, \bm{x}^{\ast}) =\{\bm{y} \ |\ \langle \bm{x}^{\ast}, \bm{y} - \bm{x} \rangle \geq 0 \}\,,
$$
$$
  H^{-}(\bm{x}, \bm{x}^{\ast}) =\{\bm{y} \ |\ \langle \bm{x}^{\ast}, \bm{y} - \bm{x} \rangle \leq 0 \}\,,
$$
and the hyperplane
$$
  H(\bm{x}, \bm{x}^{\ast}) = \{\bm{y} \ |\ \langle \bm{x}^{\ast}, \bm{y} - \bm{x} \rangle = 0 \}\,,
$$
\rsp.
Noting that in the degenerate case $ \bm{x}^{\ast} =\bm{0} $, we just have
$$
  H^{+}(\bm{x}, \bm{0}) = H^{-}(\bm{x}, \bm{0}) = H(\bm{x}, \bm{0}) = \Rn \,.
$$
To give prominence to the main theory and reduce the length of the text, we will leave all proofs of the lemmas in Appendix.

%%%%%%%%%%%%%%%%%%%%%%%%%%%%%%%%%%%%%%%%%%%%%%%%%%%%%%%%
\section{The Characterization for Convex Regularization}
%%%%%%%%%%%%%%%%%%%%%%%%%%%%%%%%%%%%%%%%%%%%%%%%%%%%%%%%

%%%%%%%%%%%%%%%%%%%%%%%%%%%%%%%%%%%%%%%%%%%%%%%%%%%%%%%%%%%%%%%%%%%%%%%%%%%%%%%%%%%%%%%%%%%%%%
\ \ \ \ \ \ {\bf Definition 3.1.}
Let $\bm{x}\in \mathrm{dom}\,f$. The set
$$
  U^{+}_{f}(\bm{x}) := \bigcup_{ \bm{x}^{\ast}\in \partial f(\bm{x})} H^{+}(\bm{x}, \bm{x}^{\ast})
$$
is called the \textbf{upper region} of $f$ at $\bm{x}$.

\vst

%%%%%%%%%%%%%%%%%%%%%%%%%%%%%%%%%%%%%%%%%%%%%%%%%%%%%%%%%%%%%%%%%%%%%%%%%%%%%%%%%%%%%%%%%%%%%%
{\bf Definition 3.2.}
Let $\,\bm{\rho}: [0,1]\ra \mathrm{dom}\, f $ be a path and $\bm{x}\in \mathrm{dom}\,f$.
The set
$$
  V^{+}_{f}(\bm{\rho}) := \bigcap_{t\in [0,1]} U^{+}_{f}(\bm{\rho(t)})
$$
is called the \textbf{ultimate region} of $f$ with respect to $\bm{\rho}$.

\vst

%%%%%%%%%%%%%%%%%%%%%%%%%%%%%%%%%%%%%%%%%%%%%%%%%%%%%%%%%%%%%%%%%%%%%%%%%%%%%%%%%%%%%%%%%%%%%%
{\bf Definition 3.3.}
A mapping $\bm{\rho}: [0,1]\ra \Rn $ is called a \textbf{searching path} with respect to $f$, if it satisfies the following conditions:

(i) $\bm{\rho} $ is continuous;

(ii) $ \bm{\rho}(t)\in \mathrm{dom}\, f $ for all $t\in [0,1]$.
%
%(iii) $f\left( \bm{\rho}(0) \right) > f \left( \bm{\rho}(1) \right) $.

\vst

%%%%%%%%%%%%%%%%%%%%%%%%%%%%%%%%%%%%%%%%%%%%%%%%%%%%%%%%%%%%%%%%%%%%%%%%%%%%%%%%%%%%%%%%%%%%%%
{\bf Remark.}
In Definition 3.3, condition (i) means the term "path" is used the same as in topology;
condition (ii) implies that $ f $ takes finite values on whole path.
%condition (iii) guarantee that the path $ \bm{\rho} $ ends at a better approximation to the minimizer of $f$ than its initial point.

\vst

%%%%%%%%%%%%%%%%%%%%%%%%%%%%%%%%%%%%%%%%%%%%%%%%%%%%%%%%%%%%%%%%%%%%%%%%%%%%%%%%%%%%%%%%%%%%%%
{\bf Definition 3.4.}
Let $\,\bm{\rho}: [0,1]\ra \mathrm{dom}\, f $ be a searching path. If
there exist a positive-valued function $\lambda: [0, \infty) \ra (0,\infty)$ and a finite-valued convex function $\psi: \Rn \ra \real $ such that
\begin{equation}\label{2301}
  \bm{\rho}(t) = \underset{\bm{x}}{\arg\min} \lbc f(\bm{x}) +\lambda(t)\psi(\bm{x}) \rbc \,.
\end{equation}
then the path $\bm{\rho}$ is said
to be \textbf{via regularization} with the \textbf{penalty} $\psi$ and the \textbf{tuner} $ \lambda $.

\vst

%%%%%%%%%%%%%%%%%%%%%%%%%%%%%%%%%%%%%%%%%%%%%%%%%%%%%%%%%%%%%%%%%%%%%%%%%%%%%%%%%%%%%%%%%%%%%%
{\bf Lemma 3.1.}
Let $ \psi $ be a finite-valued convex function and let $\lambda: [0, \infty) \ra (0,\infty)$ be a positive-valued function.
Suppose that $\,\bm{\rho}: [0,1]\ra \mathrm{dom}\, f $ is a searching path such that \eqref{2301} holds.
Then for each $ t \in [0,1] $ there exist
$$
  f'(\bm{\rho}(t))\in \partial f(\bm{\rho}(t)) \qaq \psi'(\bm{\rho}(t))\in \partial \psi(\bm{\rho}(t))
$$
such that
\begin{equation}\label{2305}
  f'(\bm{\rho}(t)) = - \lambda(t) \psi'(\bm{\rho}(t))\,.
\end{equation}
%%%%%%%%%%%%%%%%%%%%%%%%%%%%%%%%%%%%%%%%%%%%%%%%%%%%%%%%%%%%%%%%%%%%%%%%%%%%%%%%%%%%%%%%%%%%%%

\vst

%%%%%%%%%%%%%%%%%%%%%%%%%%%%%%%%%%%%%%%%%%%%%%%%%%%%%%%%%%%%%%%%%%%%%%%%%%%%%%%%%%%%%%%%%%%%%%
{\bf Lemma 3.2.}
Let $ \bm{a},\ \bm{b}\in \Rn\,,\ \theta\in (0,1)\,,\ \bm{c} = (1-\theta)\bm{a} + \theta\,\bm{b}$ and $ \delta \in (0, \theta)\,.$
Then
\begin{equation*}\label{32523}
   \bm{b} + \mathrm{cone}\left( \bm{b} - B(\bm{a}, \delta)\right)
   \subset \bm{c} + \mathrm{cone}\left( \bm{c} - B(\bm{a}, \delta)\right) \,.
\end{equation*}

\vs

Now we propose the geometric characterization of a searching path via
realization as the following theorem.

%%%%%%%%%%%%%%%%%%%%%%%%%%%%%%%%%%%%%%%%%%%%%%%%%%%%%%%%%%%%%%%%%%%%%%%%%%%%%%%%%%%%%%%%%%%%%%
{\bf Theorem 3.1.}
Let $ \psi $ be a finite-valued convex function and let $\lambda: [0, \infty) \ra (0,\infty)$ be a positive-valued function.
Suppose that $\,\bm{\rho}: [0,1]\ra \mathrm{dom}\, f $ is a searching path via regularization with penalty $\psi$ and tuner $ \lambda $.
Then the following properties hold:

(i) If $ f(\bm{\rho}(t_1)) = f(\bm{\rho}(t_2)) $,
then
\begin{equation}\label{32505}
  \psi(\bm{\rho}(t_1)) = \psi(\bm{\rho}(t_2))\,.
\end{equation}
Moreover, there is a hyperplane $ H $ such that,
for each $ t $ satisfying $f(\bm{\rho}(t)) = f(\bm{\rho}(t_1)) $,
there is a subgradients
$ f'(\bm{\rho}(t)) \in \partial f(\bm{\rho}(t)) $
such that
\begin{equation}\label{32501}
  H(\bm{\rho}(t), f'(\bm{\rho}(t))) = H \,.
\end{equation}

(ii) If $ f(\bm{\rho}(t_1)) > f(\bm{\rho}(t_2)) $, then
\begin{equation}\label{32503}
  \psi(\bm{\rho}(t_1)) < \psi(\bm{\rho}(t_2))
\end{equation}
and
\begin{equation*}\label{32502}
  \bm{\rho}(t_1)\in \mathrm{int}\,U^{+}_{f}(\bm{\rho}(t_2)) \,.
\end{equation*}

(iii) %Moreover, if $$
$ \mathrm{int}\,V^{+}_{f}(\bm{\rho}) \not= \emptyset \,. $
%
%(iv) $\bm{\rho}$ is decreasing.
%%%%%%%%%%%%%%%%%%%%%%%%%%%%%%%%%%%%%%%%%%%%%%%%%%%%%%%%%%%%%%%%%%%%%%%%%%%%%%%%%%%%%%%%%%%%%%

\vs

\zm
%%%%%%%%%%%%%%%%%%%%%%%%%%%%%%%%%%%%%%%%%%%%%%%%%%%%%%%%%%%%%%%%%%%%%%%%%%%%%%%%%%%%%%%%%%%%%%
(i).
If \eqref{32505} is not true, say,
$ \psi(\bm{\rho}_1) < \psi(\bm{\rho}_2)\,, $
then
$$
  f(\bm{\rho}_1)+ \lambda(t_2)\psi(\bm{\rho}_1) <
  f(\bm{\rho}_2)+ \lambda(t_2)\psi(\bm{\rho}_2)\,.
$$
This contradicts to the fact that $\bm{\rho}_2 $ is a minimizer of
$ f+ \lambda(t_2)\psi $. Thus \eqref{32505} must hold.

Now we show that
$$
  \mathrm{ri}\left( \mathrm{lev}_{ \leq f(\bm{\rho}_1)} f\right) = \mathrm{lev}_{ < f(\bm{\rho}_1)} f
$$
and
$$
  \mathrm{ri}\left( \mathrm{lev}_{ \leq \psi(\bm{\rho}_1)} \psi \right) = \mathrm{lev}_{ < \psi(\bm{\rho}_1)} \psi
$$
have no point in common.
If it is not true, then there is a point $ \bm{\xi} $
such that $ f(\bm{\xi}) < f(\bm{\rho}_1) $ and $  \psi(\bm{\xi}) < \psi(\bm{\rho}_1) \,. $
Thus
$$
  f(\bm{\xi})+ \lambda(t_1)\psi(\bm{\xi}) <
  f(\bm{\rho}_1)+ \lambda(t_1)\psi(\bm{\rho}_1)\,,
$$
which contradicts to the optimal property of $ \bm{\rho}_1 \,.$

Because both $ \mathrm{lev}_{\leq f(\bm{\rho}_1)} f $ and $ \mathrm{lev}_{ \leq \psi(\bm{\rho}_1)} \psi $  are non-empty convex sets, there exists a vector $\bm{u}\in \Rn $ such that the hyperplane $ H( \bm{\rho}_1, \bm{u} ) $ separates
$ \mathrm{lev}_{\leq f(\bm{\rho}_1)} f $
and $ \mathrm{lev}_{\leq \psi(\bm{\rho}_1)} \psi $.
For each $ t $ satisfying $f(\bm{\rho}(t)) = f(\bm{\rho}(t_1)) $,
$ \bm{\rho}(t)\in H( \bm{\rho}_1, \bm{u} ) $.
%Of course we can choose $ \bm{u} $ such that
%$ \mathrm{lev}_{ \leq f(\bm{\rho}_1)} f\subset H^{-}(\bm{\rho}_1, \bm{u}) \,.$
%%for all $\ \bm{x}\in \mathrm{lev}_{ \leq f(\bm{\rho}_1)} f\,, $
%%
According to Corollary 23.7.1 of \cite{rockafellar1970},
there must exist $ \lambda_t \in \real\setminus \{0\} $ such that
$$
  \bm{u} = \lambda_t f'(\bm{\rho}(t))
$$
for some $ f'(\bm{\rho}(t))\in \partial  f(\bm{\rho}(t)) $.

So we have
$$
  H ( \bm{\rho}_1, \bm{u} ) = H (\bm{\rho}(t), f'(\bm{\rho}(t)))\,,
$$
which is exactly \eqref{32501}.
%%%%%%%%%%%%%%%%%%%%%%%%%%%%%%%%%%%%%%%%%%%%%%%%%%%%%%%%%%%%%%%%%%%%%%%%%%%%%%%%%%%%%%%%%%%%%%

(ii)
If \eqref{32503} is not true, i.e.
$ \psi(\bm{\rho}_2) \leq \psi(\bm{\rho}_1)\,, $
then
$$
  f(\bm{\rho}_2)+ \lambda(t_1)\psi(\bm{\rho}_2) <
  f(\bm{\rho}_1)+ \lambda(t_1)\psi(\bm{\rho}_1)\,,
$$
a contradiction to that $\bm{\rho}_1 $ is a minimizer of
$ f+ \lambda(t_1)\psi $, and \eqref{32503} follows.

By Lemma 3.1, there exist
$ f'(\bm{\rho}_{2})\in \partial f(\bm{\rho}_{2}) $
and $ \psi'(\bm{\rho}_{2})\in \partial \psi(\bm{\rho}_{2}) $
such that
$$
  f'(\bm{\rho}_{2}) = - \lambda(t_{2})\psi'(\bm{\rho}_2)\,.
$$
Then from \eqref{32503} we have
\begin{eqnarray*}
  \bm{\rho}_1
  & \in    & \mathrm{lev}_{< \psi (\bm{\rho}_{2})} \psi \
   \subset \ \mathrm{int}\, H^{-}(\bm{\rho}_2, \psi'(\bm{\rho}_{2}))
  \ = \ \mathrm{int}\, H^{+}(\bm{\rho}_2, -\psi'(\bm{\rho}_{2}))   \\
  & =       & \mathrm{int}\, H^{+}(\bm{\rho}_2, f'(\bm{\rho}_{2}))
  \ \subset \ \mathrm{int}\, U^{+}_{f}(\bm{\rho}_2)  \,.
\end{eqnarray*}
%%%%%%%%%%%%%%%%%%%%%%%%%%%%%%%%%%%%%%%%%%%%%%%%%%%%%%%%%%%%%%%%%%%%%%%%%%%%%%%%%%%%%%%%%%%%%%

(iii).
If $ f\circ \bm{\rho}$ is constant over $[0,1]$,
so is $ \psi\circ \bm{\rho} $.
By (i) we can find a hyperplane $ H \left( \bm{\rho}(0), \bm{\rho}(0)^{\ast} \right) $
such that for each $t\in [0,1]$ there is a
$ \bm{\rho}(t)^{\ast} $ satisfying
$$
  H \left( \bm{\rho}(t), \bm{\rho}(t)^{\ast} \right) = H \left( \bm{\rho}(0), \bm{\rho}(0)^{\ast} \right)\,.
$$
Clearly, it yields
$$
  H^{+} \left( \bm{\rho}(t), \bm{\rho}(t)^{\ast} \right) = H^{+}  \left( \bm{\rho}(0), \bm{\rho}(0)^{\ast} \right)\,.
$$
and thus
$$
  H^{+}\left( \bm{\rho}(0), \bm{\rho}(0)^{\ast} \right) \subset U^{+}_{f}(\bm{\rho}(t))\,,
$$
for all $t\in [0,1]$.
Thus
$$
  H^{+}\left( \bm{\rho}(0), \bm{\rho}(0)^{\ast} \right)
  \subset \bigcap_{t\in [0,1]} U^{+}_{f}(\bm{\rho(t)})
  = V^{+}_{f}(\bm{\rho})\,,
$$
i.e. the interior of $ V^{+}_{f}(\bm{\rho}) $ contains an open half-space and, of course, is non-empty.

Now we suppose that $ f\circ \bm{\rho}$ is not constant.
Let
$$
\bm{a}\in \underset{\bm{x}}{\arg\min}\,f(\bm{x})
\qaq \bm{b}\in \underset{\bm{x}}{\arg\min}\,\psi(\bm{x})\,.
$$
Since $ f\circ \bm{\rho} $ is continuous on $[0,1]$,
we can choose $ t_0, t_1 \in[0,1] $ such that
$$
  f\circ \bm{\rho}(t_0) = \underset{ t\in[0,1] }{ \sup }\, f\circ \bm{\rho}(t)\,, \ \ \
  f\circ \bm{\rho}(t_1) = \underset{ t\in[0,1] }{ \inf }\, f\circ \bm{\rho}(t) \,,\ \
$$
\rsp. From
$$
  f(\bm{\rho}(t_0)) + \lambda(t_0) \psi(\bm{\rho}(t_0)) \leq  f(\bm{b}) + \lambda(t_0) \psi(\bm{b})
$$
and $ \psi(\bm{\rho}(t_0)) \geq \psi(\bm{b}) $,
we can see that
\begin{equation*}\label{32508}
   f\circ \bm{\rho}(t_0) \leq f(\bm{b})\,.
\end{equation*}
Similarly, we have
\begin{equation*}\label{32509}
   f\circ \bm{\rho}(t_1) \geq   f(\bm{a})\,.
\end{equation*}

Now we can choose a number $ \tau\in [0,1] $ such that
$$
   f(\bm{\rho}(t_0)) > f(\bm{\rho}(\tau)) >  f(\bm{\rho}(t_1))
   \qaq
   \psi(\bm{\rho}(t_0)) < \psi(\bm{\rho}(\tau)) < \psi(\bm{\rho}(t_1)) \,.
$$
Set
$$
  \varepsilon_{f} = \min \left( f\circ \bm{\rho}(t_0) - f\circ \bm{\rho}(\tau),
  \ f \circ \bm{\rho}(\tau) - f \circ \bm{\rho}(t_1) \right) > 0 \,,
$$
$$
  \varepsilon_{\psi} = \min \left( \psi \circ \bm{\rho}(t_1) - \psi\circ \bm{\rho}(\tau),
  \ \psi \circ \bm{\rho}(\tau) - \psi \circ \bm{\rho}(t_0) \right) > 0 \,,
$$
and
$$
  \varepsilon = \min \left( \varepsilon_{f}, \varepsilon_{\psi} \right) \,.
$$
Then
\begin{equation*}\label{32510}
   f(\bm{a}) + \varepsilon \, < f(\bm{\rho}(\tau)) \, < f(\bm{b}) - \varepsilon \,,
\end{equation*}
\begin{equation*}\label{32513}
   \psi(\bm{a}) - \varepsilon \, > \psi(\bm{\rho}(\tau)) \, > \psi(\bm{b}) + \varepsilon \,.
\end{equation*}
Now we can choose a small number $ \delta > 0 $ such that
$$
  f(\bm{x}) < f(\bm{a}) + \varepsilon \,, \ \ \text{for any}\ \bm{x} \in B(\bm{a}, \delta)\,,
$$
and
$$
  \psi(\bm{y}) < \psi(\bm{b}) + \varepsilon \,, \ \ \text{for any}\ \bm{x} \in B(\bm{b}, \delta).
$$
Define
\begin{equation*}\label{32511}
   C = \bm{b} + \mathrm{cone}\,\left( \bm{b} - B(\bm{a}, \delta) \right)
\end{equation*}
and
$$
  K = C \cap B( \bm{b}, \varepsilon ) \,.
$$

For any $ t\in [0,1] $,
if
$$
  f \circ \bm{\rho}(t) \leq f(\bm{\rho}(\tau)) \,,
$$
then
$$
  \psi \circ \bm{\rho}(t) \geq \psi(\bm{\rho}(\tau)) > \psi(\bm{b}) + \varepsilon \,.
$$
According to Lemma 3.1, there is a $ f'(\bm{\rho}(t)) \in \partial f(\bm{\rho}(t)) $
such that
$$
  -\frac{f'(\bm{\rho}(t))}{\lambda(t)} \in \partial \psi(\bm{\rho}(t))
$$
and
\begin{eqnarray}\label{32531}
   K
     & \subset & B( \bm{b}, \varepsilon )
     \subset \ \mathrm{lev}_{< \psi(\bm{b}) + \varepsilon }\,\psi \
     \subset \ \mathrm{lev}_{< \psi(\bm{\rho}(\tau))}\,\psi \
     \subset \ \mathrm{lev}_{< \psi(\bm{\rho}(t))}\,\psi
     \nonumber \\
     & \subset & H^{-}\left( \bm{\rho}(t), -\frac{f'(\bm{\rho}(t))}{\lambda(t)} \right)
     = H^{+}\left( \bm{\rho}(t), f'(\bm{\rho}(t)) \right) \
     \subset \ U^{+}_{f}\left( \bm{\rho}(t) \right) \,.
\end{eqnarray}

On the other hand,
if
$$
  f \circ \bm{\rho}(t) > f(\bm{\rho}(\tau))  \,,
$$
then for any $ f'(\bm{\rho}(t)) \in \partial f(\bm{\rho}(t)) $
we have
\begin{eqnarray*}
   B( \bm{a}, \varepsilon )
   &\subset &  \mathrm{lev}_{< f(\bm{a}) + \varepsilon }\,f \
    \subset \ \mathrm{lev}_{< f(\bm{\rho}(\tau))}\,f \
    \subset \ \mathrm{lev}_{< f(\bm{\rho}(t))}\,f
    \nonumber \\
   & \subset &  H^{-}\left( \bm{\rho}(t), f'(\bm{\rho}(t)) \right)　
     \,.
\end{eqnarray*}
Since $\bm{a}$ and $\bm{b}$ are in the opposing half-spaces associated with
the hyperplane $ H\left( \bm{\rho}(t), f'(\bm{\rho}(t)) \right) \,,$
we can find a point
$$
 \bm{c}\in \overline{\bm{a}\bm{b}}\cap H\left( \bm{\rho}(t), f'(\bm{\rho}(t)) \right) \,.
$$
By Lemma 3.2 we have
\begin{equation*}\label{32523}
   C = \bm{b} + \mathrm{cone}\left( \bm{b} - B(\bm{a}, \delta)\right)
   \subset \bm{c} + \mathrm{cone}\left( \bm{c} - B(\bm{a}, \delta)\right) \,.
\end{equation*}
Noting that
$$
  \bm{c} + \mathrm{cone}\left( \bm{c} - B(\bm{a}, \delta)\right) \subset H^{+}\left( \bm{\rho}(t), f'(\bm{\rho}(t)) \right)
$$
we have
\begin{equation}\label{32529}
  K \subset C
   \subset  H^{+}\left( \bm{\rho}(t), f'(\bm{\rho}(t)) \right) \ \subset \ U^{+}_{f}\left( \bm{\rho}(t) \right) \,.
\end{equation}
Combining \eqref{32529} with \eqref{32531}, we can see
$$
  K \ \subset \ U^{+}_{f}\left( \bm{\rho}(t) \right)
$$
for all $ t\in [0,1]$. Obviously, $ \mathrm{int}\,K $ is non-empty, which completes the proof.
\eop

%%%%%%%%%%%%%%%%%%%%%%%%%%%%%%%%%%%%%%%%%%%%%%%%%%%%
\section{The Existence of Convex Regularization}
%%%%%%%%%%%%%%%%%%%%%%%%%%%%%%%%%%%%%%%%%%%%%%%%%%%%

\ \ \ \ \ \ {\bf Definition 4.1.}
Let $\Sigma$ be an $n-$manifold in $\ \real^{n+1}\,\ \ \Sigma \subset D\times \real,\ \ D\in \Rn\,.$
We define the \textbf{steepness} of $\Sigma$ as
$$
  \mathrm{Stp}(\Sigma) := \sup\left\{ \left. \frac{|v - u|}{\| \bm{v}- \bm{u} \|}\ \right| \ (\bm{u},u),(\bm{v},v)\in \Sigma,\ \bm{v}\not= \bm{u} \right\}\,.
$$

{\bf Definition 4.2.}
Let $ f: D \ra \real $ be a continuous function, where $D\in \Rn$.
We define the \textbf{steepness} of $f$ as the steepness of its graph, $G = \{(\bm{x}, f(\bm{x}))\ |\ \bm{x}\in D\}$, i.e.
$$
  \mathrm{Stp}(f):=\mathrm{Stp}(\mathrm{graph} f)\,.
$$

{\bf Definition 4.3.}
Let $T$ be a compact and convex subset of a hyperplane $\Pi \subset \real^{n+1}$ and $\widehat{\bm{a}} \in \real^{n+1}$.
We define the \textbf{truncated cone} generated by $T$ and $\widehat{\bm{a}} $ by
$$
  \mathrm{trunc}(\widehat{\bm{a}},\ T) = \mathrm{conv}(\{\widehat{\bm{a}}\},\ T) = \{ (1-\lambda)\widehat{\bm{a}} +\lambda \widehat{\bm{x}} \ |\ \widehat{\bm{x}}\in T,\ \lambda\in [0,1] \}\,.
$$
%The point $ \widehat{\bm{a}}$, the set $T$, and the relative boundary $\mathrm{rbd}(T)$ are called the \textbf{apex}, \textbf{base}, and \textbf{directrix} of the truncated cone, respectively.
The set
$$
  \{ (1-\lambda)\widehat{\bm{a}} +\lambda \widehat{\bm{x}} \ |\ \widehat{\bm{x}}\in \mathrm{rbd}(T),\ \lambda\in [0,1] \}
$$
is called the \textbf{lateral} of the truncated cone.
%Each of the line segments between the directrix and apex is called a \textbf{generatrix} of the lateral.

\vst
%%%%%%%%%%%%%%%%%%%%%%%%%%%%%%%%%%%%%%%%%%%%%%%%%%%%%%%%%%%%%%%%%%%%%%%%%%%%%%%%%%%%
{\bf Lemma 4.1.}
Suppose that $D $ is a compact convex set in $\real^{n}$,  $\bm{a} \in \mathrm{int}(D)$, $ p, h\in \real$.
Let $S$ be the lateral of $\mathrm{trunc}(\widehat{\bm{a}}, \mathrm{rbd}(T))$, where
$$
 \widehat{\bm{a}} = (\bm{a}, p) \,, \ \ \ T = \{ (\bm{x}, h) \ |\ \bm{x}\in D\} \,.
$$
Then
$$
  \mathrm{Stp}(S) < \pinf \,.
$$

\vs

%%%%%%%%%%%%%%%%%%%%%%%%%%%%%%%%%%%%%%%%%%%%%%%%%%%%%%%%%%%%%%%%%%%
{\bf Definition 4.4.}
Let $K\subset \real^{n+1}$ be a closed convex set and $ c \in \real$ .
The set
$$
  \mathrm{sect}_{c}\,K := K \cap \{(\bm{x}, c)\ |\ \bm{x}\in \Rn \}.
$$
is called the \textbf{section} of $K$ at level $c$.
The relative boundary
$$
  \mathrm{rbd}\left(\mathrm{sect}_{c}\,K \right) =
  \left( \mathrm{sect}_{c}\,K \right) \setminus \mathrm{int}\,K
$$
is called the \textbf{transversal} of $K$ at level $c$.

\vs

%%%%%%%%%%%%%%%%%%%%%%%%%%%%%%%%%%%%%%%%%%%%%%%%%%%%%%%%%%%%%%%%%%%%%%%%%%%%%%%%%%%%
{\bf Definition 4.5.}
Let $D_1,\ \ D_2$ be compact convex sets in $\real^{n}$ and
$
  T_1 = \{(\bm{x}, h_1)\ |\ \bm{x}\in D_1\}\,,\ \ T_2 = \{(\bm{x}, h_2)\ |\ \bm{x}\in D_2\}\,,
$
where $h_1, h_2 \in\real $ with $ h_1 < h_2 $.
The set $ \mathrm{conv}\left( T_1, T_2 \right) $
is called the \textbf{frustum} generated by $T_1$ and $T_2$,
while $T_1$ and $T_2$ are called the \textbf{bottom} and \textbf{top} of frustum, respectively.
The set
$$
  \left( \mathrm{bdry}\, F \right) \setminus ( \mathrm{ri}\, T_1 \bigcup \mathrm{ri}\, T_2 )
$$
is called the \textbf{lateral} of the frustum.
The frustum $ \mathrm{conv}\left( T_1, T_2 \right) $ is said to be \textbf{top-heavy}, if
$ D_1\subset \mathrm{Int}\, D_2 \,.$

\vs

%%%%%%%%%%%%%%%%%%%%%%%%%%%%%%%%%%%%%%%%%%%%%%%%%%%%%%%%%%%%%%%%%%%%%%%%%%%%%%%%%%%%
{\bf Lemma 4.2.}
Let $D_1,\ \ D_2, \ \  h_1,\ \ h_2, \ \ T_1$, and $T_2$
be specified as in Definition 4.5.
Let $ P_{\Rn}: \real^{n+1}\ra \Rn $
be the projection operator such that $ P_{\Rn}(\bm{x}, y) = \bm{x}$ for all $ \bm{x}\in \real^{n}$ and $y\in \real $.
If the frustum $ F = \mathrm{conv}\left( T_1, T_2 \right) $ is top-heavy, then
the following hold:

(i) $ \mathrm{sect}_{h_1}\,F = T_1 $, $ \mathrm{sect}_{h_2}\,F = T_2 $;
%$ \mathrm{proj}_{\Pi_2} F = T_2 $, where $ \Pi_2 $ is the hyperplane $ \{ (\bm{x}, h_2) \ |\ \bm{x}\in \Rn \} $ and $ \mathrm{proj}_{\Pi_2}$ is the projection to $\Pi_2$;

(ii) If $ h_1 \leq c_1 < c_2 \leq h_2 $, then  $ P_{\Rn}\left( \mathrm{sect}_{c_1}\,F \right) \subset \mathrm{int} \left( P_{\Rn}\left( \mathrm{sect}_{c_2}\,F \right) \right) $;

(iii) For any $\bm{x} \in D_2\setminus \mathrm{int}\, D_1 $,
there is a unique $ y\in \real $ such that $ (\bm{x}, y)$ at the lateral of $F$.

\vs

%%%%%%%%%%%%%%%%%%%%%%%%%%%%%%%%%%%%%%%%%%%%%%%%%%%%%%%%%%%%%%%%%%%%%%%%%%%%%%%%%%%%
{\bf Definition 4.6.}
Let $D_1,\ \ D_2, \ \  h_1,\ \ h_2, \ \ T_1$, and $T_2$
be specified as in Definition 4.5 with $ D_1\subset \mathrm{int}\,D_2 $.
Let $ F = \mathrm{conv}\left( T_1, T_2 \right) $ be the top-heavy frustum generated by $T_1$ and $T_2$.
Denote by $ \Pi_{1} $ and $ \Pi_{2} $ the hyperplanes
$\{(\bm{x}, h_1)\ |\ \bm{x}\in \Rn\}$ and $\{(\bm{x}, h_2)\ |\ \bm{x}\in \Rn\}$, respectively.
Let $T \subset \Pi_{2} $ be a closed convex set and $\widehat{\bm{a}}=(\bm{a}, a)\in \real^{n+1}$ such that $ \bm{a}\in \mathrm{int}\,D_1 $ and $a<h_1$.
A truncated cone $ \mathrm{trunc}(\widehat{\bm{a}},\ T) $ is called a \textbf{upper envelope} of $F$,
if
$$
  F\subset \mathrm{trunc}(\widehat{\bm{a}},\ T) \qaq T = T_2;
$$
while
$ \mathrm{trunc}(\widehat{\bm{a}},\ T) $ is called a \textbf{lower envelope} of $F$,
if
$$
  F\subset \mathrm{trunc}(\widehat{\bm{a}},\ T) ,\ \ \ T_2\subset T \qaq \mathrm{trunc}(\widehat{\bm{a}},\ T)\bigcap \Pi_{1} = T_1 \,.
$$

\vs

%%%%%%%%%%%%%%%%%%%%%%%%%%%%%%%%%%%%%%%%%%%%%%%%%%%%%%%%%%%%%%%%%%%%%%%%%%%%%%%%%%%%
{\bf Lemma 4.3.}
Suppose that $D_1,\ \ D_2, \ \ D,\ \  h_1,\ \ h_2,\ \  \Pi_1,\ \ \Pi_2, \ \ T_1,\ \ T_2,\ \ T $
and $F$ are as in Definition 4.6.
Let $ S $ be the lateral of $F$.
Then the following hold:

(i)  If there is an upper envelop $\mathrm{trunc}(\widehat{\bm{a}},\ T_2)$ of $ F $,
then  $ \mathrm{Stp}(S) \leq \mathrm{Stp}(S_{a}) $,
where $ S_{a} $ is the lateral of $ \mathrm{trunc}(\widehat{\bm{a}},\ T_2) $.

(ii) If there is a lower envelope $\mathrm{trunc}(\widehat{\bm{b}},\ T)$ of $ F $,
then  $ \mathrm{Stp}(S) \geq \mathrm{Stp}(S_{b}) $,
where $ S_{b} $ is the lateral of $ \mathrm{trunc}(\widehat{\bm{b}},\ T) $.

%%%%%%%%%%%%%%%%%%%%%%%%%%%%%%%%%%%%%%%%%%%%%%%%%%%%%%%%%%%%%%%%%%%%%%%%%%%%%%%%%%%%

\vst

%%%%%%%%%%%%%%%%%%%%%%%%%%%%%%%%%%%%%%%%%%%%%%%%%%%%%%%%%%%%%%%%%%%
{\bf Lemma 4.4.}
Suppose that $ D_1 , \ D_2 \subset \Rn $ are bounded convex sets, $D_1\subset \mathrm{int}(D_2)\,.$
Let $ \psi: D_1\ra\real $ be a convex function, such that

(a)  $\ \psi|_{\mathrm{bdry}\,D_1} $ is constant;

(b)  $\ \mathrm{Stp}(\psi) <\pinf\,.$

\noindent Then
there is a convex function $\widetilde{\psi}: D_2\ra\real$ such that

(i) $\ \widetilde{\psi}|_{D_1} = f\,; $

(ii) $\ \widetilde{\psi}|_{\mathrm{bdry}\,D_2} $ is constant;

(iii) $\ \mathrm{Stp}(\widetilde{\psi}) <\pinf\,.$

\vst

\vst

%%%%%%%%%%%%%%%%%%%%%%%%%%%%%%%%%%%%%%%%%%%%%%%%%%%%%%%%%%%%%%%%%%%%%%%%%%%%%%%%%%%%
{\bf Lemma 4.5. }
Suppose that $ D $ is a compact convex set in $\Rn \,.$ Let $ \psi: D \ra\real $ be a convex function, such that

  1)  $\ \psi|_{\partial D} = c $ is a constant;

  2)  $\ \mathrm{Stp}(\psi) < \pinf $.

\noindent Then there is a convex function $\widetilde{\psi}:\Rn\ra \real$ such that
$\widetilde{\psi}|_{D} = \psi \,.$

\vs

%%%%%%%%%%%%%%%%%%%%%%%%%%%%%%%%%%%%%%%%%%%%%%%%%%%%%%%%%%%%%%%%%%%%%%%%%%%%%%%%%%%%

\vs

%%%%%%%%%%%%%%%%%%%%%%%%%%%%%%%%%%%%%%%%%%%%%%%%%%%%%%%%%%%%%%%%%%%%%%%%%%%%%%%%%%%%
{\bf Theorem 4.1. }
Let $m$ be a positive integer and $ \bm{x}_{i}\in \mathrm{dom}\,f \,,\ \ i = 0, 1, \cdots, m $.
Suppose there are subgradients
$ \bm{x}_{i}^{\ast} \in \partial\,f(\bm{x}_{i}) \,,\ \ i = 0, 1, \cdots, m\,, $
satisfying that the following conditions:

(i) If $ f(\bm{x}_i) > f(\bm{x}_j) $, then $ \bm{x}_i\in \mathrm{int}\,H^{+}(\bm{x}_j, \bm{x}_{j}^{\ast}) $;

(ii) If $ f(\bm{x}_i) = f(\bm{x}_j) $, then $ H^{+}(\bm{x}_i, \bm{x}_{i}^{\ast}) = H^{+}(\bm{x}_j, \bm{x}_{j}^{\ast}) $;

(iii) The set $ \mathrm{int}\bigcap_{ i= 1}^{n} H^{+}(\bm{x}_i, \bm{x}_{i}^{\ast}) $ is nonempty.

\noindent Then there are a convex function $\ \psi $ and a positive numbers $\lambda_{i},\ \ i = 1,2,\cdots, n$, such that
\begin{equation*}\label{2721}
  \bm{x}_{i} = \underset{\bm{x}}{\arg\min} \lbc f(\bm{x}) +\lambda_i \psi(\bm{x}) \rbc \ ,\ \ \ i = 1,\cdots, m .
\end{equation*}

\vs

%%%%%%%%%%%%%%%%%%%%%%%%%%%%%%%%%%%%%%%%%%%%%%%%%%%%%%%%%%%%%%%%%%%%%%%
\zm
\walog, we suppose %that the sequence $ \left( f(\bm{x}_{0}), f(\bm{x}_{1}), \cdots, f(\bm{x}_{m}) \right) $ is decreasing, i.e.
\begin{equation}\label{32201}
  f(\bm{x}_{0}) =\cdots = f(\bm{x}_{i_{1}-1})> f(\bm{x}_{i_{1}}) = \cdots = f(\bm{x}_{i_{2}-1})> \cdots >  f(\bm{x}_{i_{k}}) = \cdots = f(\bm{x}_{m})\,,
\end{equation}
where $ 0 < i_{1} < i_{2} < \cdots <i_{k} \leq m = i_{k+1}-1 $ for some positive integer $ k \leq m $.

According to condition (iii), we can choose a point
$\bm{a} \in \mathrm{int} \bigcap_{i=0}^{m} H^{+}(\bm{x}_i, \bm{x}_{i}^{\ast}) $
and a small number $\ \varepsilon > 0 $
such that the closed ball of center $\bm{a}$ and radius $ \varepsilon $
\begin{equation*}\label{32303}
  \overline{B(\bm{a}, \varepsilon)}\subset \mathrm{int} \bigcap_{i=0}^{m} H^{+}(\bm{x}_i, \bm{x}_{i}^{\ast}) \,.
\end{equation*}

Noting that condition (ii) combined with \eqref{32201} implies
\begin{equation*}\label{32301}
  H(\bm{x}_{0}, \bm{x}_{0}^{\ast}) = \cdots = H(\bm{x}_{i_{1}-1}, \bm{x}_{i_{1}-1}^{\ast})\,,
\end{equation*}
we can see
\begin{equation*}\label{32304}
  \bm{x}_{j} \in H(\bm{x}_{0}, \bm{x}_{0}^{\ast}),\quad \text{for}\ j = 0, 1, \cdots, i_{1}-1 \,.
\end{equation*}
From condition (ii) and \eqref{32201} we deduce
\begin{equation}\label{32305}
  \bm{x}_{j} \in \mathrm{int}\bigcap_{h = i_{1}}^{m} H^{+}(\bm{x}_{h}, \bm{x}_{h}^{\ast}),\quad \text{for}\ j = 0, 1, \cdots, i_{1}-1 \,.
\end{equation}
Thus we can choose two balls
$ B_{0,s} $ and $ B_{0,e} $
such that
$ \overline{B}_{0,s} $ is tangent to $ H(\bm{x}_0, \bm{x}_{0}^{\ast})$ at $ \bm{x}_0$,
$ \overline{B}_{0,e} $ is tangent to $ H(\bm{x}_{0}, \bm{x}_{0}^{\ast})$ at $ \bm{x}_{i_{1}-1}$,
and
\begin{equation}\label{32306}
  \overline{B_{0,s}}\cup \overline{B_{0,e}} \subset \bigcap_{h = i_{1}}^{m} H^{+}(\bm{x}_{h}, \bm{x}_{h}^{\ast}) \,,
\end{equation}
where "s" and "e" refer to "start" and "end", \rsp.

Then we denote
\begin{equation*}\label{32302}
  C_0 = \mathrm{conv}\left( \overline{B}_{0,s}\cup \overline{B}_{0,e}\cup \{ \bm{x}_{1}, \cdots, \bm{x}_{i_{1}-2} \} \right) \,.
\end{equation*}
From \eqref{32305} and \eqref{32306} we can see
\begin{equation*}\label{32307}
  C_0 \subset \mathrm{int}\bigcap_{h = i_{1}}^{m} H^{+}(\bm{x}_{h}, \bm{x}_{h}^{\ast}) \,.
\end{equation*}
Denote
\begin{equation}\label{2722}
  K_{0} = \mathrm{conv}\left( \overline{B(\bm{a} \cup \varepsilon)}, C_0 \right)
  \,.
\end{equation}
Clearly, $ K_{0} $ is compact, convex, and satisfying
\begin{equation*}\label{32308}
  K_{0}\subset \mathrm{int}\bigcap_{h = i_{1}}^{m} H^{+}(\bm{x}_{h}, \bm{x}_{h}^{\ast}) \,.
\end{equation*}

In general, for $j = 0,1, \cdots, k-1$, if a compact convex set
\begin{equation}\label{32309}
  K_{j} \subset \mathrm{int} \bigcap_{h = i_{(j+1)}}^{m} H^{+}(\bm{x}_{h}, \bm{x}_{h}^{\ast})
\end{equation}
is determined,
then we define
\begin{equation}\label{32101}
  E_{j} = \{ \bm{x}\ |\ \mathrm{dist}\left( \bm{x}, K_{j}\right) \leq \frac{d_j}{2} \} \,,
\end{equation}
where
$$
  d_{j} = \mathrm{dist}\left( K_{j},\ \mathrm{bdry} \bigcap_{h = i_{(j+1)}}^{m} H^{+}(\bm{x}_{h}, \bm{x}_{h}^{\ast}) \right) \,.
$$
Since $ K_{j} $ is compact and, according to \eqref{32309},
we can assert $ d_{j} > 0 $.
Clearly, $E_{j}$ is also compact and convex.
From \eqref{32309} and \eqref{32101} we can see that
\begin{equation*}\label{32310}
  E_{j} \subset \mathrm{int} \bigcap_{h = i_{(j+1)}}^{m} H^{+}(\bm{x}_{h}, \bm{x}_{h}^{\ast}) \,.
\end{equation*}

Now we select two balls
$ B_{j+1, s} $ and $ B_{j+1, e} $ such that
$ \overline{B}_{j+1, s} $ is tangent to $ H(\bm{x}_{i_{(j+1)}}, \bm{x}_{i_{(j+1)}}^{\ast})$ at $ \bm{x}_{i_{(j+1)}}$,
$ \overline{B}_{j+1, e} $ is tangent to $ H(\bm{x}_{i_{(j+2)}-1}, \bm{x}_{i_{(j+2)}-1}^{\ast})$ at $ \bm{x}_{i_{(j+2)}-1}$,
and
\begin{equation*}\label{32366}
  \overline{B_{j+1, s}}\cup \overline{B_{j+1, e}} \subset \bigcap_{h = i_{(j+1)}}^{m} H^{+}(\bm{x}_{h}, \bm{x}_{h}^{\ast}) \,.
\end{equation*}
Then we denote
\begin{equation*}\label{32302}
  C_{j+1} = \mathrm{conv}\left( \overline{B}_{j+1, s}\cup \overline{B}_{j+1, e}\cup \{ \bm{x}_{i_{(j+1)}+1}, \cdots, \bm{x}_{i_{(j+1)}-2} \} \right) \,.
\end{equation*}

We set
\begin{equation*}\label{32102}
  K_{j+1} = \mathrm{conv}\left( E_{j} \cup C_{j+1} \right)\,.
\end{equation*}
Then, $ K_{j+1} $ is compact, convex, and satisfying
\begin{equation*}\label{32311}
  K_{j+1}\subset \mathrm{int} \bigcap_{h = i_{(j+2)}}^{m} H^{+}(\bm{x}_i, \bm{x}_i^{\ast}) \,,
\end{equation*}
where we define $\ \bigcap_{h = i_{(k+1)}}^{m} H^{+}(\bm{x}_i, \bm{x}_i^{\ast}) = \Rn\,. $

In this way we recursively construct a series of compact convex sets
$ K_{0}, K_{1}, \cdots, K_{k}$
with the following properties:

a) $\ K_{j} \subset \mathrm{int}\,K_{j+1} $, for $ j =0, 1, \cdots, k-1\, ;$

b) $\ K_{j} $ is tangent to $\ H(\bm{x}_{i_{j}}, \bm{x}_{i_{j}}^{\ast})$ at $ \bm{x}_{i_{j}}, \cdots, \bm{x}_{i_{(j+1)}-1}$, for $ j =0, 1, \cdots, k\,.$

Noting that $ \bm{a}\in \mathrm{int}\,K_{0}$ by \eqref{2722}, we define a function $ \psi_{0}: K_{0}\ra \real $
as follows:
For any $ \bm{x}\in K_{0} $, by applying the convexity of $K_{0}$, we can find a point $ \bm{y}\in \mathrm{bdry}\,K_{0} $ such that $ \bm{x} = (1-\lambda)\bm{a} +\lambda \bm{y} $ for some $\lambda\in [0,1]$ and then simply define
\begin{equation*}\label{32103}
  \psi_{0}(\bm{x}) = \lambda \,.
\end{equation*}
Note that $\psi_{0}(\bm{a}) = 0 $ and $ \psi_{0}|_{\mathrm{bdry}K_{0}} = 1 $.
Since $ \mathrm{epi}\,\psi_{0} $ is the intersection of the convex cone
$$
  \{ \left. \left( \bm{a} + t(\bm{x}-\bm{a}), t \right) \in\real^{n+1} \ \right| \ \bm{x}\in K_{0},\ t \geq 0 \}
$$
and the convex column
$$
  \{ \left. \left( \bm{x}, z \right) \in\real^{n+1} \ \right|\ \bm{x}\in K_{0},\ z \in \real \}\,,
$$
$ \mathrm{epi}\,\psi_{0} $ is convex and so is $\psi_{0}$.
Obviously, the graph of $ \psi_{0} $ is the lateral of
$ \mathrm{trunc}\left( (\bm{a}, 0), T_0 \right) $, where
$ T_0 = \{(\bm{x}, 1)\ |\ \bm{x}\in K_{0}\} $.
By Lemma 4.1 we have $ \mathrm{Stp}(\psi_{0}) < \pinf \,. $

For $j = 1, 2, \cdots, k$, applying Lemma 4.4 recursively,
we can obtain convex function
$ \psi_{j}: K_{j}\ra \real $
such that
$$
  \psi_{j}|_{K_{j-1}} = \psi_{j-1}\,,\ \ \ \psi_{j}|_{\mathrm{bdry}K_{j}} = c_{j}\,,\qaq \mathrm{Stp}(\psi_{j}) <\pinf \,,
$$
where $c_{j}$ is a constant.

By Lemma 4.5, $\psi_{m} : K_{k}\ra \real $ can be extended to a convex function $\psi: \Rn \ra \real \,.$

The property b) means that, at each $ \bm{x}_{i} $, $ i = 0, 1, \cdots, m$, the isosurface of $f$,
$$
   \{\bm{x}\ |\ f(\bm{x}) = f(\bm{x}_{i})\} \,,
$$
is tangent to the isosurface of $\psi$,
$$
  \{\bm{x}\ |\ \psi(\bm{x}) = \psi(\bm{x}_{i})\}\,.
$$
Since
$$
  \arg\min\ \psi \in \mathrm{int}\,H^{+}(\bm{x}_{i}, \bm{x}_{i}^{\ast}) \,,
$$
there is a $f'(\bm{x}_{i}) \in \partial f(\bm{x}_{i})$ that is opposite to some $\psi'(\bm{x}_{i}) \in \partial \psi(\bm{x}_{i})$, for each $ i \in \{0, 1, \cdots, m \}\,.$
Thus there are real numbers $ \lambda_{i}$ such that
\begin{equation*}%\label{2730}
  f'(\bm{x}_{i}) = - \lambda_i \psi'(\bm{x}_{i}), \ \ i = 0, 1, \cdots, m\,.
\end{equation*}
Thus $\bm{x}_{i}$ is a critical point of $\ f(\bm{x}) +\lambda_i \psi(\bm{x})$ for $ i=0, 1, \cdots, m\,$, which completes the proof.
\eop

%%%%%%%%%%%%%%%%%%%%%%%%%%%%%%%%%%%%%%%%%%%%%%%%%%%%%%%%%%%%%%%%%%%%%%%%%%%%%%%%%%%%
\section{The Approximation by Convex Regularization}
%%%%%%%%%%%%%%%%%%%%%%%%%%%%%%%%%%%%%%%%%%%%%%%%%%%%

{\bf Theorem 5.1. }
Let $\,\bm{\rho}: [0,1]\ra \mathrm{dom}\, f $ be a searching path. Suppose that there exist a positive-valued function $\lambda: [0, \infty) \ra (0,\infty)$ and a finite-valued convex function $ \psi $ such that
the following conditions are satisfied:

i) If $ f(\bm{\rho}(t_2)) < f(\bm{\rho}(t_1)) $, then
$ \bm{\rho}(t_1)\in \mathrm{int}\left( U^{+}_{f}(\bm{\rho}(t_2)) \right) $;

ii) $ \mathrm{int}\left(  V^{+}_{f}(\bm{\rho}) \right) \not= \emptyset $.

Then for any $ \varepsilon > 0 $, there are a convex function $\ \psi $ and a positive-valued function $\lambda: [0,1]\ra (0, \pinf)$, such that
for each $ t\in [0,1] $, there is a minimizer of $ f(\cdot)+\lambda(t) \psi(\cdot) $
\begin{equation}\label{3001}
  \bm{x}(t) \in \underset{\bm{x}}{\arg\min} \lbc f(\bm{x}) +\lambda(t) \psi(\bm{x}) \rbc
\end{equation}
such that
\begin{equation}\label{3002}
  \| \bm{\rho}(t) - \bm{x}(t) \| < \varepsilon\,.
\end{equation}

\zm
Since $\bm{\rho}$ is uniformly continuous on $[0,1]$, we can choose a natural number $ m $ and a partition of $[0,1]$, $0 = t_0 < t_1 < \cdots < t_m = 1$,
such that
\begin{equation*}\label{3003}
  \| \bm{\rho}(t) - \bm{\rho}(t_j) \| < \varepsilon\,, \ \ \text{ for } t\in [t_{j-1}, t_j],\ \ j=1,2,\cdots,m\,.
\end{equation*}
By Theorem 3.1, we have a convex function $ \psi $ and a positive numbers $\lambda_{j},\ \ j = 1,2,\cdots, m$, such that
\begin{equation*}\label{3004}
  \bm{x}_{j} \in \underset{\bm{x}}{\arg\min} \lbc f(\bm{x}) +\lambda_i \psi(\bm{x}) \rbc \ ,\ \ \ j = 1,\cdots, m .
\end{equation*}
We define $\lambda: [0,1]\ra (0, \pinf)$ and $\bm{x}: [0,1]\ra (0, \pinf)$
by
\begin{equation*}\label{3005}
  \left\{
    \begin{array}{c}
      \lambda(t) = \lambda_{j}\,, \\
      \bm{x}(t) = \bm{x}_{j}\,,
    \end{array}
  \right.
  \ \ \  t \in [t_{j-1}, t_j]\,, \ \ j = 1, 2, \cdots, m \,.
\end{equation*}
Then both \eqref{3001} and \eqref{3002} hold.
\eop

\section{Conclusion}

We show that there is an intimate connection between penalization and early stopping.
In fact, it is almost a necessary and sufficient condition under which a search path of a convex optimization problem can be represented by a penalization course. In this way one can study the statistical features of an iterative algorithm by exploring the correspondent penalization function, which is easier to be handled quantitatively than an algorithmic course.

\bibliographystyle{unsrt}
 % \bibliographystyle{plainnat}
 %\bibliographystyle{elsarticle-num}
 %\bibliographystyle{elsarticle-harv}
 %\bibliographystyle{elsarticle-num-names}
 %\biboptions{square,numbers,sort&compress}
%%%%%%%%%%%%%%%%%%%%%%%%%%%%%%%%%%%%%%%%%%%%%%%%%%%%%%%%%%%%%%%%%%%%%%%%%%%%%%%%%%%%%%%%%%%%%%%%%
\bibliography{mybib}
%%%%%%%%%%%%%%%%%%%%%%%%%%%%%%%%%%%%%%%%%%%%%%%%%%%%%%%%%%%%%%%%%%%%%%%%%%%%%%%%%%%%%%%%%%%%%%%%

\newpage

%%%%%%%%%%%%%%%%%%%%%%%%%%%%%%%%%%%%%%%%%%%%%%%%%%%%%%%%%%%%%%%%%%%%%%%%%%%%%%%%%%%%%%%%%%%%%%%%%
\begin{center}
{\bf  Appendix}
\end{center}
%%%%%%%%%%%%%%%%%%%%%%%%%%%%%%%%%%%%%%%%%%%%%%%%%%%%%%%%%%%%%%%%%%%%%%%%%%%%%%%%%%%%%%%%%%%%%%%%%
\vst

{\bf  Proof of Lemma 3.1. }
From \eqref{2301} we have
\begin{equation*}%\label{2304}
  \bm{0} \in \partial \left( f(\bm{\rho}(t)) +\lambda(t)\psi(\bm{\rho}(t)) \right)\,, \ \  t\in [0,1]\,.
\end{equation*}
According to Theorem 23.8, Rockfella70, p.223,
$$
  \partial \left( f(\bm{\rho}(t)) +\lambda(t)\psi(\bm{\rho}(t)) \right) = \partial f(\bm{\rho}(t)) + \lambda(t) \partial \psi(\bm{\rho}(t))\,.
$$
Thus there exist
$
  f'(\bm{\rho}(t))\in \partial f(\bm{\rho}(t)) \qaq \psi'(\bm{\rho}(t))\in \partial \psi(\bm{\rho}(t))
$
such that
$$
  f'(\bm{\rho}(t)) + \lambda(t) \psi'(\bm{\rho}(t)) = 0\,,
$$
which is equivalent to \eqref{2305}.
\eop

\vst

{\bf  Proof of Lemma 3.2. }
We first assume that $ \bm{a} = \bm{0} \,.$
Denote
\begin{eqnarray*}
  K_{\bm{b}} &=&  \bm{b} + \mathrm{cone}\left( \bm{b} - B(\bm{0}, \delta)\right)\,,  \\
  K_{\bm{c}} &=&  \bm{c} + \mathrm{cone}\left( \bm{c} - B(\bm{0}, \delta)\right) \,.
\end{eqnarray*}
For any $\bm{x} \in K_{\bm{b}}\,,$ there exist $\bm{u}\in B(\bm{0}, \delta) $ and  $\lambda \geq 0 $  such that
\begin{eqnarray*}
  \bm{x}
      &=& \bm{b} + \lambda \left( \bm{b} - \bm{u} \right) \\
      &=& \theta\,\bm{b} + \frac{1+\lambda - \theta}{\theta} \left( \theta \bm{b} - \frac{\lambda\theta}{1+\lambda-\theta} \bm{u} \right) \\
      &=& \bm{c} + \widetilde{\lambda }\left( \bm{c} - \widetilde{ \bm{u}} \right) \,,
\end{eqnarray*}
where
$$
  \widetilde{\lambda } =  \frac{1+\lambda - \theta}{\theta} > 0  \qaq
  \widetilde{ \bm{u}} = \frac{\lambda\theta}{1+\lambda-\theta} \bm{u} \,.
$$
Since $ 0 < \theta < 1 $, we have
$ (1+\lambda)\theta < 1+\lambda $ or $ \lambda \theta < 1+\lambda - \theta $.
Thus
$$
  0 < \frac{\lambda\theta}{1+\lambda-\theta} < 1 \,.
$$
Hence $ \| \widetilde{ \bm{u}} \| < \| \bm{u} \| < \delta $
and $ \widetilde{ \bm{u}}\in B(\bm{0}, \delta) \,.$
Thus $ \bm{x} \in K_{\bm{c}} \,. $

In general case, noting that
\begin{eqnarray*}
  \bm{b} + \mathrm{cone}\left( \bm{b} - B(\bm{a}, \delta)\right)
  &=& \bm{a} + (\bm{b}- \bm{a}) + \mathrm{cone}\left( \bm{b} - \bm{a} + B(\bm{0}, \delta)\right) \\
  &=& \bm{a} + \widetilde{\bm{a}} + \mathrm{cone}\left( \widetilde{\bm{b}} - B(\bm{0}, \delta)\right) \,,  \\
  \bm{c} + \mathrm{cone}\left( \bm{c} - B(\bm{a}, \delta)\right)
  &=& \bm{a} + (\bm{c}- \bm{a}) + \mathrm{cone}\left( \bm{c} - \bm{a} + B(\bm{0}, \delta)\right) \\
  &=& \bm{a} + \widetilde{\bm{c}} + \mathrm{cone}\left( \widetilde{\bm{c}} - B(\bm{0}, \delta)\right) \,,
\end{eqnarray*}
where
$$
  \widetilde{\bm{b}} = \bm{b} - \bm{a}\,,\ \ \  \widetilde{\bm{c}} = \bm{c} - \bm{a}\,,
$$
we can obtain
$$
  \bm{b} + \mathrm{cone}\left( \bm{b} - B(\bm{a}, \delta)\right)
  \ =\ \bm{a}+ K_{\widetilde{\bm{b}}}
  \ \subset\ \bm{a}+ K_{\widetilde{\bm{c}}}
  \ =\ \bm{c} + \mathrm{cone}\left( \bm{c} - B(\bm{a}, \delta)\right)
$$
from $ K_{\widetilde{\bm{b}}} \subset K_{\widetilde{\bm{c}}} $.
\eop

\vs

{\bf  Proof of Lemma 4.1. }
Since $D $ is compact, so is $\mathrm{bd}(D)$.
Noting that $ \bm{a} \not\in \mathrm{bd}(D) $,
we have $\mathrm{dist}(\bm{a},\mathrm{rbd}(D)) >0 $ and
then
$$
  M := \underset{\bm{x}\in \mathrm{rbd}(D)}{ \sup } \frac{|p-h|}{\| \bm{x}-\bm{a} \|}=\frac{|p-h|}{\mathrm{dist}(\bm{a},\mathrm{rbd}(D))} < \pinf\,.
$$
For any two points $ \widehat{\bm{y}}_1 ,\ \widehat{\bm{y}}_2\in S $, there are
$ \bm{x}_1,\ \bm{x}_2 \in \mathrm{bd}(D) $ and $ \lambda_1,\ \lambda_2 \in [0,1] $ such that
$$
  \widehat{\bm{y}}_1 = (1-\lambda_1)\widehat{\bm{x}}_1 + \lambda_1\widehat{\bm{a}},\ \ \ \widehat{\bm{y}}_1 = (1-\lambda_2)\widehat{\bm{x}}_2 + \lambda_2 \widehat{\bm{a}},
$$
where $  \widehat{\bm{x}}_1 = (\bm{x}_1, h),\ \ \widehat{\bm{x}}_2 = (\bm{x}_2, h) \in \mathrm{rbd}(T) $.
Thus
$
  \widehat{\bm{y}}_1 = (\bm{y}_1, k_1),\ \ \widehat{\bm{y}}_2 = (\bm{y}_2, k_2)\,,
$
where
\begin{eqnarray*}
  \bm{y}_1 &=& (1-\lambda_1)\bm{x}_1 + \lambda_1 \bm{a},\ \ \ \bm{y}_2 = (1-\lambda_2)\bm{x}_2 + \lambda_2 \bm{a}\,, \\
  k_1 &=& (1-\lambda_1)h + \lambda_1 p,\ \ \ k_2 = (1-\lambda_2)h + \lambda_2 p \,.
\end{eqnarray*}
Then we have
\begin{equation*}
  \frac{|k_2 - k_1|}{\| \bm{y}_2 - \bm{y}_1 \|} =
  \frac{|\left( (1-\lambda_1)h + \lambda_1 p \right) - \left( (1-\lambda_2)h + \lambda_2 p \right)|}{\|\left( (1-\lambda_1)\bm{x}_1 + \lambda_1 \bm{a} \right) - \left( (1-\lambda_2)\bm{x}_2 + \lambda_2 \bm{a} \right) \|}
  =  \frac{|p-h|}{\bm{\xi}-\bm{a}}\,,
\end{equation*}
where
$$
  \bm{\xi}=\frac{1-\lambda_1}{\lambda_2-\lambda_1}\bm{x}_1 +  \frac{\lambda_2-1}{\lambda_2-\lambda_1}\bm{x}_2 \,.
$$
We denote
$$
  \alpha = \frac{1-\lambda_1}{\lambda_2-\lambda_1} , \ \ \ \beta = \frac{\lambda_2 -1}{\lambda_2-\lambda_1} \,.
$$
Then
$$
  \bm{\xi} = \alpha\bm{x}_1 +\beta\bm{x}_2 \,.
$$
Note that $ \alpha + \beta=1 $,
which means that $\bm{\xi} $ is at the line passing through $\bm{x}_1$ and $\bm{x}_2$.
On the other hand, $\lambda_2 -1\leq 0$, which yields
$$
 \alpha \beta = (1-\lambda_1)(\lambda_2 -1)/(\lambda_2-\lambda_1)^2 \leq 0 \,,
$$
so that  $\bm{\xi}\not\in \overline{\bm{x}_1\bm{x}_2} $.
Thus $ \bm{\xi}\not\in \mathrm{int}D $.
Then $ \| \bm{\xi} -\bm{a} \| \geq \mathrm{dist}(\bm{a},\mathrm{rbd}(D)) $, so that
$$
  \frac{|k_2 - k_1|}{\| \bm{y}_2 - \bm{y}_1 \|} \leq \frac{|p-h|}{\mathrm{dist}(\bm{a},\mathrm{rbd}(D))} = M <\pinf \,.
$$
\eop

\vs

{\bf  Proof of Lemma 4.2. }

(i): Trivial.

(ii): We first show that $ P_{\Rn}\left( \mathrm{sect}_{c_1}\,F \right) \subset \mathrm{int}\,D_{2} $.

%$ P_{\Rn}\left( \mathrm{sect}_{c_1}\,F \right) \subset P_{\Rn}\left( \mathrm{sect}_{c_2}\,F \right) $.
When $c_1 = h_1$, it is directly from (i), (ii) and the top-heavy assumption.
So we need only to discuss the case $ h_1< c_1 < h_2 $.

If $ \bm{x}\in P_{\Rn}\left( \mathrm{sect}_{c_1}\,F \right) $, then $ (\bm{x}, c_1)\in F $.
Thus there are a pair of natural numbers $k,m$, a set of points $ \bm{x_1},\cdots, \bm{x_k}\in D_1 $,
a set of points $ \bm{x_{m+1}},\cdots, \bm{x_{k+m}}\in D_2 $, and nonnegative numbers
$ \lambda_{1}, \cdots, \lambda_{k+m}\in [0,1] $ such that $ \lambda_{1}+ \cdots+ \lambda_{k+m}=1 $ and
$$
  \lambda_{1}(\bm{x_1},h_1) + \cdots + \lambda_{k}(\bm{x_k},h_1)+\lambda_{m+1}(\bm{x_{m+1}},h_2) + \cdots + \lambda_{k}(\bm{x_{k+m}},h_2) = (\bm{x}, c_1)\,.
$$
It implies
$$
  (\lambda_{1}+\cdots+\lambda_{k})h_1 + (\lambda_{k+1}+\cdots+\lambda_{k+m})h_2 = c_1 = (\lambda_{1}+\cdots+\lambda_{k+m})c_1 \,.
$$
Noting that $ h_1 < c_1 < h_2 $ implies
$$
  \lambda_{1}+\cdots+\lambda_{k} > 0 \qaq \lambda_{k+1}+\cdots+\lambda_{k+m}>0 \,,
$$
we can write
$$
  \bm{x} = (\lambda_{1}+\cdots+\lambda_{k})\bm{\xi}_{1} +(\lambda_{k+1}+\cdots+\lambda_{k+m})\bm{\xi}_{2}\,,
$$
where
$$
  \bm{\xi}_{1} = \frac{\lambda_{1}}{\lambda_{1}+\cdots+\lambda_{k}}\bm{x}_{1}+\cdots+\frac{\lambda_{k}}{\lambda_{1}+\cdots+\lambda_{k}}\bm{x}_{k}\in D_1\subset \mathrm{int}\,D_{2}\,,
$$
$$
  \bm{\xi}_{2} = \frac{\lambda_{k+1}}{\lambda_{k+1}+\cdots+\lambda_{k+m}}\bm{x}_{k+1}+\cdots+\frac{\lambda_{k+m}}{\lambda_{k+1}+\cdots+\lambda_{k+m}}\bm{x}_{k+m}\in D_{2}\,.
$$
Since $\lambda_{1}+\cdots+\lambda_{k} > 0$ and $\bm{\xi}_{1}\in \mathrm{int}\,D_{2}\,,$ we can conclude that
$
  \bm{x} \in \mathrm{int}\,D_{2}\,.
$

Now we turn to prove that
\begin{equation}\label{31401}
  \bm{x} \in \mathrm{int}\, P_{\Rn}\left( \mathrm{sect}_{c_2}\,F \right) \,.
\end{equation}

If $c_2 = h_2$, in this case we have $P_{\Rn}\left( \mathrm{sect}_{c_2}\,F \right) = D_2 $ and the inclusion relation has already been established.
So we only consider the case $ c_2 < h_2 $.

Because $ \bm{x} $ is in the interior of $D_2$,
there is ball $B(\bm{x}, r)$ of center $\bm{x}$ and radius $ r > 0 $ such that
$ B(\bm{x}; r)\subset \mathrm{int}\,D_2 $.
Thus the set $ \widehat{B} = \{(\bm{x}, h_2)\ |\ \bm{x}\in  B(\bm{a}; r)\}\subset T_2 $.
Then the truncated cone
$$
  \mathrm{trunc}\left( (\bm{x}, c_1), \widehat{B} \right) \subset F \,.
$$
Noting that $ c_1 < c_2 < h_2 $, we can deduce that
$$
  (\bm{x}, c_2) \in \mathrm{int}\,\mathrm{trunc}\left( (\bm{a}, h_1), B(\bm{a}; r) \right) \subset \mathrm{int}\,F \,.
$$
It is equivalent to \eqref{31401}.

(iii).
If it is not true,
then there exist some $\bm{a}\in D_2\setminus \mathrm{int}\, D_1 $,
$c_1$ and $c_2 \in \real $ with $c_1 < c_2$, such that
\begin{equation*}\label{31301}
    (\bm{a}, c_1), (\bm{a}, c_2)\in \left( \mathrm{bdry}\, F \right) \setminus ( \mathrm{ri}\, T_1 \bigcup \mathrm{ri}\, T_2 ) \,,
\end{equation*}
which implies
$$
   (\bm{a}, c_1) \in \mathrm{ri}\left( \mathrm{sect}_{c_1} F \right)
   \qaq
   (\bm{a}, c_2) \in \mathrm{ri}\left( \mathrm{sect}_{c_2} F \right) \,.
$$
Thus
$$
  \bm{a}\in \mathrm{bdry}\left( P_{\Rn}\left( \mathrm{sect}_{c_1}\,F \right) \right)
  \qaq
  \bm{a}\in \mathrm{bdry}\left( P_{\Rn}\left( \mathrm{sect}_{c_2}\,F \right) \right)\,.
$$
It contradicts to (iii), which asserts
$ \bm{a}\in \mathrm{int}\left( P_{\Rn}\left( \mathrm{sect}_{c_2}\,F \right) \right) $ for
$ \bm{a}\in \mathrm{bdry}\left( P_{\Rn}\left( \mathrm{sect}_{c_1}\,F \right) \right) $.
\eop

\vs

{\bf  Proof of Lemma 4.3. }
(i). Let
$$
 \widehat{\bm{u}}= (\bm{u}, u), \ \  \widehat{\bm{v}}= (\bm{v}, v) \in S\,,
$$
where $ \bm{u}, \bm{v}\in \Rn,\ \ u,v \in \real $. Without any loss of generality we assume $ u < v< h_2 $.
Let
$ \widehat{\bm{w}} = (\bm{w}, h_2) $ be the intersection of $\mathrm{ray}(\widehat{\bm{u}}, \widehat{\bm{v}})$ and $\Pi_2$, i.e.
$$
  \{ \widehat{\bm{w}} \} = \mathrm{ray}(\widehat{\bm{u}}, \widehat{\bm{v}}) \cap \Pi_2.
$$

We first show that
\begin{equation}\label{2703}
  \bm{w}\not\in \mathrm{int}\,D_2 \ .
\end{equation}
In fact, if $\bm{w} \in \mathrm{Int}(D_2)$, then we can choose a small number
$ \varepsilon > 0 $, such that the set
$$
   B = \{(\bm{x}, h_2)\ |\ \| \bm{x}- \bm{w}\| < \varepsilon \}\subset T_2 \ .
$$
Since $ \widehat{\bm{u}}\in \mathrm{Conv}(T_1, T_2) $, we have
\begin{equation}\label{2704}
  \mathrm{trunc}(\widehat{\bm{u}},B) \subset F\ .
\end{equation}
Remember that $v < h_2$, which combining with \eqref{2704} implies that $\widehat{\bm{v}}\in \mathrm{int}\,F$.
This contradicts to the fact that $ \widehat{\bm{v}}\in S \subset \mathrm{bdry}\,F $, and \eqref{2703} follows.

Since $ \bm{u}\in \mathrm{int}\,D_2 $, there exists a point $ \bm{w}_{\ast} $ such that $ \bm{w}_{\ast} \in \mathrm{bdry}\, D_2 $
and $ \bm{w}_{\ast}\in \overline{\bm{u}\bm{w}} $.
Thus
$ \| \bm{w}_{\ast} - \bm{u} \| \leq \| \bm{w} - \bm{u} \|\,.$
On the other hand, because
$ F $ is contained in $\mathrm{trunc}(\widehat{\bm{a}},\ T_2)$,
we have
$
 \mathrm{sect}_{u} F \subset\mathrm{sect}_{u}\left( \mathrm{trunc}(\widehat{\bm{a}},\ T_2) \right) \,
$
and then
$$
 P_{\Rn}\left( \mathrm{sect}_{u} F \right) \subset P_{\Rn}\left(  \mathrm{sect}_{u}\left( \mathrm{trunc}(\widehat{\bm{a}},\ T_2) \right) \right) \,.
$$
Thus there exists a point $ \bm{u}_{\ast} $ such that
$$
 \bm{u}_{\ast} \in \overline{\bm{u}\bm{w}_{\ast}}\, \cap\, \mathrm{bdry}\, P_{\Rn}\left(  \mathrm{sect}_{u}\left( \mathrm{trunc}(\widehat{\bm{a}},\ T_2) \right) \right) \,.
$$
Noting that by Lemma 4.2 we have
$$
 P_{\Rn}\left(  \mathrm{sect}_{u}\left( \mathrm{trunc}(\widehat{\bm{a}},\ T_2) \right) \right) \subset \mathrm{int}\,
 P_{\Rn}\left(  \mathrm{sect}_{h_2}\left( \mathrm{trunc}(\widehat{\bm{a}},\ T_2) \right) \right)
 = \mathrm{int}\,D_2
$$
and thus
$ \bm{u}_{\ast} \in \mathrm{int}\,D_2 $, in particular, $ \bm{u}_{\ast} \not= \bm{w}_{\ast}\in \mathrm{bdry}\,D_2 $.

Note that $ (\bm{u}_{\ast}, u),\ (\bm{w}_{\ast}, h_2)\in S_{a} $.
Since $ (\bm{u}, u),\ (\bm{u}_{\ast}, u),\ (\bm{v}, v) $ and $ (\bm{w}, h_2)$ are colinear,
we have
$$
  \frac{|v-u|}{\| \bm{v} - \bm{u} \|} = \frac{|h_2 - u|}{\| \bm{w} - \bm{u} \|}
   \leq \frac{|h_2 - u|}{\| \bm{w}_{\ast} - \bm{u}_{\ast} \|}
   \leq \mathrm{Stp}(S_a)\,,
$$
and then
\begin{equation*}
  \mathrm{Stp}(S) = \underset{(\bm{x},u),(\bm{y},v)\in S,\ u < v}{\sup}\frac{|v - u|}{\| \bm{v}- \bm{u} \|} \leq \mathrm{Stp}(S_a)\,.
\end{equation*}

(ii). Similarly to (i).
\eop

\vs

{\bf  Proof of Lemma 4.4. }
Let $\psi|_{\mathrm{bdry}\,D_1} \equiv c_1 \in \real .$
Denote
$$
  G_1 = \{(\bm{x}, \psi(\bm{x}))\ |\ \bm{x}\in D_1\}\ ,\ \ \ T_1 = \{(\bm{x}, c_1)\ |\ \bm{x}\in D_1\}\ ,\ \ \ F_1 = \mathrm{conv}\,G_1,\,.
$$
Select arbitrarily a point $ \bm{a}\in \mathrm{int}\,D_1 $. Since $ \mathrm{bdry}\,D_1 $ is compact, we can find a point $ \bm{x}_1\in \mathrm{bdry}\,D_1 $ such that
$$
  \| \bm{x}_1 - \bm{a} \| = \underset{\bm{x}\in \mathrm{bdry}\,D_1}{\sup}\mathrm{dist}(\bm{x}, \bm{a}).
$$
Then we choose a number $ p\in\real $ such that $ p < c_1 + \mathrm{Stp}(\psi) \cdot \| \bm{x}_1 - \bm{a} \| \,,$ which equivalent to
$$
  \frac{c_1 - p}{\| \bm{x}_1 - \bm{a} \|} > \mathrm{Stp}(\psi) \,.
$$
Denote
$$
  \widehat{\bm{a}} = (\bm{a}, p) \qaq K_1 = \mathrm{trunc}(\widehat{\bm{a}},\ T_1)\,.
$$
We take five steps to prove the lemma.

\textbf{Step 1.}
We show that
\begin{equation}\label{31404}
  F_1 \subset K_1\,.
\end{equation}

In fact, if it is not true, then there is a point at $G_1$, the graph of $\psi$, and in the exterior of $ K_1 $, i.e. there is some $ \bm{b}\in D_1 $ such that
\begin{equation}\label{31402}
  (\bm{b}, \psi(\bm{b}))\not\in K_1\,.
\end{equation}
This also provides $ \bm{b}\not=\bm{a} $ and there is a point
$ \bm{b}_{\ast}\in \mathrm{bdry}\left( P_{\Rn}\left( \mathrm{sect}_{f(\bm{b}) }K_1 \right) \right) $
such that $ \bm{b}_{\ast}\in \mathrm{ray}(\bm{a},\bm{b}) $.
Noting that \eqref{31402} implies $ \bm{b}\not\in P_{\Rn}\left( \mathrm{sect}_{f(\bm{b}) }K_1 \right) $, thus $ \bm{b}_{\ast}\in \overline{\bm{a}\bm{b}} $.
On the other hand, there is a point
$ \bm{b}^{\ast}\in \mathrm{bdry} D_1 $.
such that $ \bm{b}^{\ast}\in \mathrm{ray}(\bm{a},\bm{b}) $.
$ \bm{b}\in \overline{\bm{a}\bm{b}^{\ast}} $,
for $ \bm{b}\in \mathrm{int}\,D_1 $.

So we have
$$
  \| \bm{b}_{\ast} - \bm{a} \| < \| \bm{b} - \bm{a} \| < \| \bm{b}^{\ast} - \bm{a} \| \,.
$$
Then, for $ (\bm{a},p) $, $ (\bm{b}_{\ast}, \psi(\bm{b})) $ and $ (\bm{b}^{\ast},h_1) $ are colinear,
\begin{equation}\label{31405}
  \frac{c_1 - \psi(\bm{b})}{\bm{b}^{\ast} - \bm{b}}\leq \mathrm{Stp}(G_1) = \mathrm{Stp}(\psi)\,.
\end{equation}
On the other hand,
$$
  \mathrm{Stp}(\psi)\leq \frac{c_1 - p}{\bm{x}_1 - \bm{a}}= \frac{c_1 - \psi(\bm{b})}{\bm{b}^{\ast} - \bm{b}_{\ast}}\,,
$$
which contradicts to \eqref{31405}. Thus \eqref{31404} has to be valid.

\textbf{Step 2.}
We define the function $\widetilde{\psi}$ with (i) and (ii) satisfied.

For every $ \bm{x}\in \mathrm{bdry}\,D_1 $,
we define
$ \rho(\bm{x}) $ to be the unique point in $ (\mathrm{bdry}\,D_2)\cap\mathrm{ray}(\bm{a},\bm{x}) $.
Noting that the function $ \frac{\|\rho(\bm{x})-\bm{a}\|}{\|\bm{x}-\bm{a}\|} $ of $ \bm{x} $ is continuous on the compact set $\mathrm{bdry}\,D_1 $,
we can choose a point $ \bm{x}^{\ast}\in \mathrm{bdry}\,D_1 $ such that
$$
  \frac{\|\rho(\bm{x}^{\ast}) - \bm{a}\|}{\|\bm{x}^{\ast} - \bm{a}\|} = \underset{\bm{x}\in \mathrm{bdry}\,D_1}{\sup}\frac{\|\rho(\bm{x})-\bm{a}\|}{\|\bm{x}-\bm{a}\|}
  \geq 1 \,.
$$
Let
$$
  \lambda^{\ast} = \frac{\|\rho(\bm{x}^{\ast}) - \bm{a}\|}{\|\bm{x}^{\ast} -\bm{a}\|}\,.
$$
Define
$$
  \tau(\bm{x}) = (1-\lambda^{\ast})\bm{a} + \lambda^{\ast}\bm{x}\,, \ \ \bm{x}\in\,\mathrm{bdry}\,D_1\,,
$$
and
$$
  c_2 = (1-\lambda^{\ast})p+\lambda^{\ast}c_1 \,.
$$
Then
\begin{equation}\label{31406}
  D_2\subset \tau(D_1)\,,
\end{equation}
since $ \rho(\bm{x})\in \overline{\bm{a}\tau(\bm{x})} $ for every $\bm{x}\in \mathrm{bdry}\,D_1 $,
and the set
\begin{equation}\label{31409}
  T = \{(\tau(\bm{a}), c_2)\ |\ \bm{x}\in \,\mathrm{bdry}\,D_1\}\,.
\end{equation}
is compact and convex, since it is simply the image of a compact convex set via an affine transformation, adding a translation in the last coordinate component.

Denote
$$
  T_2 = \{(\bm{x}, c_2)\ |\ \bm{x}\in D_2\} \qaq F_2 = \mathrm{conv}(T_1, T_2).
$$
Then $F_2$ is a top-heavy truncated frustum.
According to Lemma 4.2.(iii), we can define
a function $ g : D_2\setminus \mathrm{int}\,D_1\,\ra \real$ as
$$
  g(\bm{x}) = y\,,\ \ \text{for each }(\bm{x},y)\in S\,,
$$
where $S$ is the lateral of $F_2$.
Thus we can define $\widetilde{\psi}: D_2 \ra \real $ as
\begin{equation*}\label{31408}
  \widetilde{\psi}(\bm{x})=
    \left\{
      \begin{array}{ll}
        f(\bm{x}), & \text{for $\bm{x}\in D_1$}; \\
        g(\bm{x}), & \text{for $\bm{x}\in D_2\setminus \mathrm{int}\,D_1$}.
      \end{array}
    \right.
\end{equation*}
Parts (i) and (ii) of the conclusion of this lemma can be verified straightforward.

\textbf{Step 3. }
We show that $\widetilde{\psi}$ is convex.

To do this, we need only to prove that
the epigraph of $\widetilde{\psi}$, or, equivalently, $F_1\cup F_2$ is convex.
We use reduction to absurdity again.
If it is false, then there exist some
$$
 \widehat{\bm{u}} = (\bm{u}, u)\in F_1 \qaq \widehat{\bm{v}} = (\bm{v}, v)\in F_2
$$
such that for some $ \lambda \in (0,1) $,
$
 (1-\lambda)\widehat{\bm{u}}+\lambda \widehat{\bm{v}}\not\in F_1\cup F_2\,.
$
Thus there are two intersection points $ \widehat{\bm{u}}^{\ast}  = (\bm{u}^{\ast}, u^{\ast}) $ and $\widehat{\bm{v}}^{\ast} = (\bm{v}^{\ast}, v^{\ast}) $ such that
$$
  \widehat{\bm{u}}^{\ast} ,\  \widehat{\bm{v}}^{\ast}
  \in \overline{\widehat{\bm{u}}\widehat{\bm{v}}}\bigcap \mathrm{bdry}(F_1 \cup F_2)
$$
and
$$
  \left( \mathrm{ri}\,\overline{\widehat{\bm{u}}^{\ast} \widehat{\bm{v}}^{\ast}} \right) \bigcap \mathrm{bdry}(F_1 \cup F_2) = \emptyset\,.
$$

Obviously, $ \widehat{\bm{u}}^{\ast}$ and $ \widehat{\bm{v}}^{\ast}$ cannot appear in only one of the truncated cones $F_1$ and $F_2$.
Thus we can assume that
$$
 \widehat{\bm{u}}^{\ast} \in (\mathrm{bdry}\,F_1)\setminus T_1 \qaq \widehat{\bm{v}}^{\ast} \in (\mathrm{bdry}\,F_2)\setminus (T_1\cup T_2)\,.
$$
In this case, however, $ \overline{\widehat{\bm{u}}^{\ast} \widehat{\bm{v}}^{\ast}} $ has a intersection point, say, $ \widehat{\bm{w}} $, with the
hyperplane $ \{ (\bm{x}, c_1)\ |\ \bm{x} \in \Rn \} $.
Since $ T_1 \in F_1 $, $ \widehat{\bm{w}}\not\in T_1 $.

On the other hand,
denote
$$
  \widetilde{K_1} = \mathrm{trunc}(\widehat{\bm{a}}, T)\,,
$$
where $T$ is defined in \eqref{31409}.
It is easy to see that $ K_1 \subset \widetilde{K_1} $.
Noting that \eqref{31406} implies $ T_2\subset T $,
we conclude that
$$
  F_1\cup F_2 \subset \widetilde{K_1}\,,
$$
for
$ F_1 \subset K_1 \subset \widetilde{K_1} $
, $ T_1, T_2 \subset \widetilde{K_1}$.
So the truncated cone $ \widetilde{K_1} $ is convex
and a lower envelope of the
top-heavy truncated frustum $F_2$.
since $\widetilde{K_1}$ is convex and
$ \widehat{\bm{u}}^{\ast},\ \widehat{\bm{v}}^{\ast}\in \widetilde{K_1}$,
we have
$ \widehat{\bm{w}}\in \widetilde{K_1}$.
So
$$
  \widehat{\bm{w}}\in \widetilde{K_1}\bigcap  \{ (\bm{x}, c_1)\ |\ \bm{x} \in \Rn \} = \mathrm{sect}_{c_1} \widetilde{K_1} = T_1\,,
$$
which is a contradiction to the relation $ \widehat{\bm{w}}\not\in T_1 $  we have just proved.
Therefore, we obtain convexity of $F_1 \cap F_2$.

\textbf{Step 4. }
we set up
\begin{equation}\label{31510}
  \mathrm{Stp}(S_2) \, < \, \pinf \,,
\end{equation}
where $S_2$ is the lateral of $F_2$.

We denote
$$
  \widehat{\bm{a}}_t = (\bm{a}, t)\,,\
  K_t = \mathrm{trunc}(\widehat{\bm{a}}_t, T_2)\,, \ T_t = \mathrm{sect}_{c_1}K_t\,,\ \text{and} \ D_t = P_{\Rn}T_t\,,\ \text{for } t < c_1 \,.
$$
Obviously,
$ T_t \subset \mathrm{ri}\,T_2\,,$ or,
equivalently, $D_t \subset \mathrm{int}\,D_2\,. $
Define
$$
  \delta(t):= \underset{\bm{x}\in \partial D_t}{\sup}\ \mathrm{dist}(\bm{x}, \partial D_2) > 0 \,.
$$
It is easy to see that
$$
  \delta(t)\ra 0\,, \ \ \text{as } t\ra \ninf \,.
$$
Since $ \mathrm{dist}(\partial D_1, \partial D_2)>0 $,
we can find some real number $ q < p $ such that
$$
  \delta(q) < \mathrm{dist}(\partial D_1, \partial D_2)\,,
$$
in which case
$
   D_1 \subset \mathrm{int}\,D_q\,,
$
or, equivalently,
$
  T_1 \subset \mathrm{ri}\,T_q\,.
$
Thus $K_q$ is an upper envelope of $F_2 = \mathrm{conv}(T_1, T_2)$.
From Lemma 4.3.(i), \eqref{31510} has been established.

\textbf{Step 5. }
The remaining is to prove $\mathrm{Stp}(\widetilde{\psi}) < \pinf $.
We need only to verify that
\begin{equation}\label{31501}
  \mathrm{Stp}(G_1 \cup S_2) < \pinf\,,
\end{equation}
because $ G_1 \cup S_2 $ is exactly the graph of $ \widetilde{\psi} $.

Let $ \widehat{\bm{u}} = (\bm{u},u) $ and $ \widehat{\bm{v}} = (\bm{v},v)$ be two points in
$ G_1 \cup S_2 $. Without of any generality we suppose that $ u < v $. We show that
\begin{equation}\label{31502}
  \frac{v-u}{\| \bm{v} - \bm{u} \|} \,\leq\, \max\left( \mathrm{Stp}(G_1),\,\mathrm{Stp}(S_2) \right)\,.
\end{equation}
Having this valid, we just obtain \eqref{31501} by taking supremum over all possible choices of $\widehat{\bm{u}} $ and $\widehat{\bm{v}} $ on the left-side of \eqref{31502}.
The discussion can be taken in three cases as follows.

Case 1: $ \widehat{\bm{u}}, \widehat{\bm{v}}\in G_1 $.
In this case we have simply
\begin{equation*}%\label{31513}
  \frac{v-u}{\| \bm{v} - \bm{u} \|} \,\leq\, \mathrm{Stp}(G_1)\,.
\end{equation*}

Case 2: $ \widehat{\bm{u}}, \widehat{\bm{v}}\in S_2 $.
Similar to Case 1, we have
\begin{equation*}%\label{31514}
  \frac{v-u}{\| \bm{v} - \bm{u} \|} \,\leq\, \mathrm{Stp}(S_2)\,.
\end{equation*}

Case 3: $ \widehat{\bm{u}}\in G_1 $ and $\widehat{\bm{v}}\in S_2 $.
In this case there is a unique point
$$
   \widehat{\bm{w}}=(\bm{w}, c_1) \in \overline{\widehat{\bm{u}}\widehat{\bm{v}}} \cap \Pi_1\,,
$$
where the hyperplane $\Pi_1 = \{(\bm{x}, c_1)\ |\ \bm{x}\in \Rn \}$.

If $ \widehat{\bm{w}} \in T_1 $, then from the colinearity of $ \widehat{\bm{u}} $, $\widehat{\bm{v}} $ and
$ \widehat{\bm{w}} $, we have
\begin{equation}\label{31505}
  \frac{v-u}{\| \bm{v} - \bm{u} \|} = \frac{v-c_1}{\| \bm{v} - \bm{w} \|} \,\leq\, \mathrm{Stp}(G_1)\,.
\end{equation}

If $ \widehat{\bm{w}} \not\in T_1 $,
then $\bm{w}\not\in D_1 $.
Thus we can find a point
$$
  \bm{w}_{\ast}\in \overline{\bm{u} \bm{w}}\cap \partial D_1 \,,
$$
which implies
\begin{equation}\label{31504}
  \| \bm{w}_{\ast} - \bm{u} \| < \| \bm{w} - \bm{u} \|\,.
\end{equation}
Thus from Equality \eqref{31504} and the colinearity of $ \widehat{\bm{u}} $, $\widehat{\bm{w}} $ and
$ \widehat{\bm{v}} $, we have
\begin{equation}\label{31506}
  \frac{v-u}{\| \bm{v} - \bm{u} \|} = \frac{v-c_1}{\| \bm{v} - \bm{w} \|} \,< \,
  \frac{v-c_1}{\| \bm{v} - \bm{w}_{\ast} \|} \,\leq\,
  \mathrm{Stp}(G_1)\,.
\end{equation}

At all events, either \eqref{31505} or \eqref{31506} holds and it follows \eqref{31502}.
It completes the proof.
\eop

\vs

{\bf  Proof of Lemma 4.5. }
Denote
$$
  G = \{(\bm{x}, \psi(\bm{x}))\ |\ \bm{x}\in D \}\ ,\ \ \ T = \{(\bm{x}, c)\ |\ \bm{x}\in D \}\,.
$$
Similar to the construction in the proof of Lemma 4.4, we can
choose a point $ \widehat{\bm{a}} = (\bm{a}, p) $
such that
$$
  \bm{a}\in \mathrm{int}\,D \,,\ \ \ G \subset \mathrm{trunc}(\widehat{\bm{a}}, T) \,.
$$
Define
$$
  K = \{ (1-t)\bm{a}+ t\bm{x}\ |\ \bm{x}\in T, \ t\geq 1 \}
$$
and
denote $S_K$ to be the lateral of $K$, i.e.
$$
  S_K = \{ (1-t)\bm{a}+ t\bm{x}\ |\ \bm{x}\in \partial T, \ t \geq 1 \}\,.
$$
Then
We can define $ \widetilde{\psi} $ as
\begin{equation*}\label{2728}
  \widetilde{\psi}(\bm{x}) =
  \left\{
    \begin{array}{ll}
      \psi(\bm{x}), & \text{for } \bm{x}\in D\,, \\
      y, &  \text{for } \bm{x}\in \Rn \setminus D \text{ and } (\bm{x}, y) \in S_K \,.
    \end{array}
  \right.
\end{equation*}
Then we can easily check that $ \widetilde{\psi} $ is well-defined.
$ \widetilde{\psi}|_{D} = f $ is straightforward obtained.
The convexity of $\widetilde{\psi}$ can be shown from a deduction similar to that applied in Step 3 of the proof of Lemma 4.4.
\eop

\end{document}